\documentclass[oneside, 11pt]{amsart}
\usepackage[utf8]{inputenc}

\usepackage{amsfonts}
\usepackage{amsmath}
\usepackage{amssymb}
\usepackage{amsthm}
\usepackage[a4paper]{geometry}
\usepackage{mathrsfs} 
\usepackage[dvipsnames]{xcolor}
\usepackage[colorlinks=true, linkcolor=blue, citecolor=blue]{hyperref}
\usepackage{appendix}
\usepackage{enumerate}
\usepackage{geometry}

\theoremstyle{plain}
\newtheorem{theorem}{Theorem}
\newtheorem{proposition}[theorem]{Proposition}
\newtheorem{lemma}[theorem]{Lemma}
\newtheorem{corollary}[theorem]{Corollary}
\newtheorem{definition}[theorem]{Definition}
\newtheorem{remark}[theorem]{Remark}

\theoremstyle{definition}

\numberwithin{theorem}{section}
\numberwithin{equation}{section} 

\newcommand{\set}[1]{\left\{ #1 \right\}}
\newcommand{\norm}[1]{\left \lVert  #1 \right \rVert}
\newcommand{\abs}[1]{\left\lvert #1 \right\rvert}

\newcommand{\Z}{\mathbb{Z}}

\newcommand{\R}{\mathbb{R}}
\newcommand{\N}{\mathbb{N}}

\newcommand{\calF}{\mathcal{F}}
\newcommand{\calG}{\mathcal{G}}
\newcommand{\calM}{\mathcal{M}}

\newcommand{\bbE}{\mathbb{E}}

\newcommand{\condexp}[2]{\bbE \left[ #1 \lvert #2 \right]}

\title[Dynamical Gibbs Variational Principle]{Dynamical Gibbs Variational Principles for Irreversible Interacting Particle Systems with Applications to Attractor Properties}
\author{Benedikt Jahnel}
\address{Institut f\"ur Mathematische Stochastik, Technische Universit\"at Braunschweig, Universit\"atsplatz 2,
38106 Braunschweig, Germany \& Weierstrass Institute for Applied Analysis and Stochastics\\
Mohrenstraße 39\\
10117 Berlin\\
Germany}
\email{benedikt.jahnel@tu-braunschweig.de}

\author{Jonas Köppl}
\address{Weierstrass Institute for Applied Analysis and Stochastics\\
Mohrenstraße 39\\
10117 Berlin\\
Germany}
\email{jonas.koeppl@wias-berlin.de}

\date{
    \today
}

\begin{document}
\maketitle
\begin{abstract}
  \noindent 
  We consider irreversible translation-invariant interacting particle systems on the $d$-dimensio\-nal cubic lattice with finite local state space, which admit at least one Gibbs measure as a time-stationary measure. Under some mild non-degeneracy conditions on the rates and the specification we prove that zero relative entropy loss of a translation-invariant measure implies that the measure is Gibbs w.r.t.~the same specification as the time-stationary Gibbs measure. As an application, we obtain the attractor property for irreversible interacting particle systems, which says that any weak limit point of any trajectory of translation-invariant measures is a Gibbs measure w.r.t.~the same specification as the time-stationary measure. This extends previously known results to fairly general irreversible interacting particle systems. 
\end{abstract}




\section{Introduction and finite state space analogy}
\subsection{Introduction}
Interacting particle systems are countable systems of locally interacting Markov processes and are often used as toy models for stochastic phenomena with an underlying spatial structure. An original motivation for studying such systems came from statistical mechanics. The idea was to describe and analyze stochastic models for the time evolution of systems whose equilibrium states are the classical Gibbs measures. In particular, one hoped to obtain a better understanding of the phenomenon of phase transitions.

Even though the definition of an interacting particle system often looks very simple and the major technical issues of its existence and uniqueness for sufficiently well-behaved transition rates have long been settled, it is in general surprisingly difficult to say anything non-trivial about their behavior. In most cases, explicit calculations are not feasible and one has to be content with qualitative statements and estimates. Some of the main challenges deal with the long-time behavior of the systems. The first step of proving any limit theorem is to describe the possible limit points of the time-evolved distribution $\nu_t$ as $t$ tends to infinity. As a next step, one can then try to determine the basin of attraction.

In the case of irreducible finite-state Markov processes, this question has long been answered, but for interacting particle systems this question is much trickier and in many situations a part of the difficulty is due to non-uniqueness of time-stationary distributions. In this regime, the analysis is very delicate and various techniques have been developed to study limit theorems or attractor properties. One particular technique that will play a major role in this manuscript is due to Holley~\cite{holley_free_1971} and involves using the relative entropy functional with respect to some specification $\gamma$ as a Lyapunov function for the measure-valued differential equation that describes the time evolution of the system in the space of measures.

This idea was later extended to more general but still reversible systems by Higuchi and Shiga~\cite{higuchi_results_1975}. A couple of years later Künsch~\cite{kunsch_non_1984} managed a first step towards a treatment of non-reversible systems, but only in the simple case of single-site updates and binary local state space. Moreover, Künsch only treated the case where the specification is given through a potential. Recently, Jahnel and Külske \cite{jahnel_attractor_2019} managed to extend the previous results to very general systems but came just short of verifying the attractor property for general non-reversible systems. It is therefore the main purpose of this manuscript to go beyond the existing literature and establish a dynamical Gibbs variational principle, plus the corresponding attractor property, for general non-reversible interacting particle systems. Conceptually our proof strategy is mainly inspired by the strategy in \cite{kunsch_non_1984} but extends the results for non-reversible systems to the more general setting in \cite{jahnel_attractor_2019}. 

Let us note that the relative entropy and its rate of change are very common tools for studying systems of interacting particles, and they connect probability, analysis, and geometry intricately. One particularly fruitful application of relative entropy techniques is in the context of Log-Sobolev inequalities for Markov processes. These inequalities can be used to obtain bounds on the (exponential) speed of convergence to equilibrium. However, these methods are limited to the situation where the time-stationary measure is unique, whereas our method goes beyond this case and is also applicable in the non-uniqueness regime. A pedagogical introduction to Log-Sobolev inequalities in the easier setting of Markov chains on finite state spaces can be found in \cite{diaconis_logarithmic_1996}, while a very general approach can be found in \cite[Chapter 5]{bakry_analysis_2014}. 

Another sub-area where relative entropy methods have successfully been applied is the derivation of hydrodynamic equations from microscopic models of interacting particle systems. In this context, the method is used to study the infinite particle limit, with additional rescaling of space and time, and not for long-time asymptotics. An introduction to this method can for example be found in the monograph \cite{kipnis_scaling_1999}. 
 
A quite recent application of relative entropy in a very similar setting as ours are the works \cite{chazottes_gaussian_2020} and \cite{chazottes_relative_2022} on Gaussian concentration and uniqueness criteria for Gibbs measures. One of the main ideas in these recent publications is to use that certain concentration properties – that are satisfied by high-temperature Gibbs measures –  imply that the relative entropy density with respect to $\mu$ is positive definite. 

It is also noteworthy that Holley’s method is not limited to interacting particle systems on the $d$-dimensional cubic lattice $\Z^d$, but has recently also been used to study systems on more general, even non-amenable, graphs, see \cite{shriver_free_2022}. 

One can also use a similar approach, involving the decay of relative entropy, to prove the central limit theorem. This was first observed by Linnik in \cite{linnik_information-theoretic_1959} for i.i.d.~$\R$-valued random variables whose law is absolutely continuous with respect to Lebesgue measure. Since then, the results have been extended to more general situations. 
In \cite{carlen_entropy_1991}, Carlen and Soffer employed the entropy functional as a Lyapunov function to study central limit theorems for sums of independent and dependent random variables and also applied their results to prove central limit theorems for a number of lattice spin systems. A couple of years later, the assumption of absolute continuity with respect to Lebesgue measure was lifted in \cite{artstein_solution_2004}, at least in the case of independent random variables. There,  it was shown that the entropy is always increasing along sequences of the form $(n^{-1/2}\sum_{i=1}^n X_i)_{n \in \N}$ where $(X_i)_{i \in \N}$ are i.i.d.~square-integrable random variables, not necessarily absolutely continuous with respect to the Lebesgue measure. 
This monotonicity of the entropy of the partial sums appearing in the classical central limit theorem can be interpreted as a formal analogue of the second law of thermodynamics. Roughly speaking, this law says that, as a thermally insulated system evolves towards thermodynamic equilibrium, its entropy is non-decreasing.

The rest of our paper is organized as follows. In Section~\ref{section:finite-state-space} we motivate our method of proof by considering the simple example of a continuous-time Markov chain on a finite state space. In Section~\ref{section:setting-and-results} we introduce the basic setting of infinite-volume Gibbs measures and interacting particle systems, before we then formulate and discuss our results. The proofs of these can then be found in Section~\ref{section:proofs}.
\subsection{Relative entropy loss in finite state spaces}\label{section:finite-state-space}
We now want to consider monotonicity properties of the entropy or free energy in a more dynamical setting, namely in the context of continuous-time Markov processes. There we will encounter a dynamic counterpart to the classical Gibbs variational principle as stated for example in ~\cite[Theorem 6.82]{friedli_statistical_2017}.
We again note a formal similarity with the second law of thermodynamics but will not enter a physical discussion thereof and instead refer the interested reader to the excellent references~\cite{landau_statistische_1987,kardar_statistical_2007} and the in-depth discussion of the second law of thermodynamics in \cite{lieb_physics_1999}.
To explain the general ideas in a simple setting, we will first restrict ourselves to the case of Markovian dynamics on a finite state space $E$. 

For two probability measures $\mu, \nu$ on $E$ with $\mu(x) >0$ for all $x\in E$ we define the \textit{relative entropy of $\nu$ with respect to $\mu$} by 
\begin{align*}
        h(\nu | \mu) := \sum_{x \in E}\nu(x) \log\left(\frac{\nu(x)}{\mu(x)}\right),
\end{align*}
where we use the convention that $0 \log 0 = 0$. 
Let $\mathscr{L}=(L_{x y})_{x,y\in E}$ be an irreducible generator for a continuous-time Markov process on $E$. In this situation, it is well known that there exists a unique measure $\mu$, which is time-stationary with respect to the Markov semigroup $(e^{t \mathscr{L}})_{t \geq 0}$, see e.g. \cite[Proposition 2.61]{liggett_continuous_2010}. Moreover, we have $\mu(x) > 0$ for all $x \in E$. Note that we do not assume that $\mu$ is reversible with respect to $\mathscr{L}$. For an initial distribution $\nu \in \calM_1(E)$ we denote the distribution at time $t \geq 0$ by $\nu_t$, i.e., 
\begin{align}
    \nu_t(B) = \int_{E}e^{t \mathscr{L}}\mathbf{1}_B(\omega) \ \nu(d\omega), \quad B \subset E. 
\end{align}
We will recall that the relative entropy can be used as a Lyapunov function, i.e., for all initial distributions $\nu \in \calM_1(E)$, the map
$
    t \mapsto h(\nu_t \lvert \mu) 
$
is non-increasing and only vanishes for $\nu = \mu$. For this, we analyze the \textit{relative entropy loss}, which is defined by
\begin{align}
    \label{def-relative-entropy-loss}
    g_{\mathscr{L}}(\nu | \mu) = \frac{d}{dt}\rvert_{t=0} \ h(\nu_t | \mu). 
\end{align}
By a simple calculation one sees that for $\nu \in \calM_1(E)$ the relative entropy loss can be written as 
\begin{equation}
    \label{relative-entropy-loss-via-generator}
        g_{\mathscr{L}}(\nu | \mu) = \sum_{x \in E}\left[\sum_{y \neq x}\nu(y) L_{y x}\log\left(\frac{\nu(x)}{\mu(x)}\right) - \sum_{y \neq x}\nu(x) L_{x y}\log\left(\frac{\nu(x)}{\mu(x)}\right)\right]. 
\end{equation}
\noindent 
Equipped with this representation of $g_{\mathscr{L}}(\nu | \mu)$ we see that it is non-positive. Indeed,  consider the function
\begin{align*}
    \Phi: \R \to \R,\qquad\Phi(u) 
    :=
    \begin{cases}
        u - u\log(u) - 1, \quad &\text{if } u > 0,
        \\\
        -1, &\text{if } u \leq 0, 
    \end{cases}
\end{align*}
and note that $\Phi$ is strictly concave on $[0,\infty)$ and only takes non-positive values.
Without loss of generality, we can assume that $\nu(x) > 0$ for all $x \in E$. Otherwise, we would have $g_{\mathscr{L}}(\nu | \mu) = - \infty \leq 0$. 
    Since $\mu$ is time-stationary with respect to the Markov process generated by $\mathscr{L}$, we know that for all $x \in E$,
    \begin{align}\label{proof:sign-relative-entropy-loss-finite-space-1}
        \sum_{y \neq x}\mu(y) L_{y x} = \sum_{y \neq x}\mu(x) L_{x y} = \mu(x) L_{x x}. 
    \end{align}
    This implies that 
    \begin{align}\label{finite-state-space-essential-equality}
        \sum_{x \in E }\sum_{y\neq x} \nu(x) \frac{\mu(y)}{\mu(x)}L_{y x}
        = 
        \sum_{x \in E}\frac{\nu(x)}{\mu(x)}\sum_{y \neq x} \mu(y) L_{y x} 
        = 
        \sum_{x \in E} \nu(x) L_{x x} 
        =
        \sum_{x \in E} \nu(x) \sum_{y\neq x} L_{x y}, 
    \end{align}
    and we  can use~\eqref{finite-state-space-essential-equality} to write 
    \begin{align*}
        g_{\mathscr{L}}(\nu | \mu) 
        =
        &\sum_{x \in E}
            \left[
                \sum_{y \neq x}
                \nu(y) L_{y x}\log\left(\frac{\nu(x)}{\mu(x)}\right) 
                - 
                \sum_{y \neq x}
                \nu(x) L_{x y}\log\left(\frac{\nu(x)}{\mu(x)}\right) 
            \right]
        \\\
        =
        & \sum_{x \in E}\sum_{y \neq x}\nu(x) L_{y x} \frac{\mu(y)}{\mu(x)} \Phi\left(\frac{\mu(x)}{\nu(x)}\frac{\nu(y)}{\mu(y)} \right).
    \end{align*}
    Since $\Phi$ is non-positive, this implies that the relative entropy loss is also non-positive, and therefore the relative entropy is non-increasing along trajectories. 
    This is already a nice and intuitive result in itself, but we can deduce even more information from the representation 
    \begin{align}\label{eqn:relative-entropy-loss-finite-state-spaces}
        g_{\mathscr{L}}(\nu | \mu) 
        =
        \sum_{x \in E}\sum_{y \neq x}\nu(x) L_{y x} \frac{\mu(y)}{\mu(x)} \Phi\left(\frac{\mu(x)}{\nu(x)}\frac{\nu(y)}{\mu(y)} \right). 
    \end{align}
    Because $\mu$ is time-stationary, it is clear that $g_\mathscr{L}(\mu | \mu) = 0$. But $\mu$ is also uniquely characterized as the maximizer of the relative entropy loss functional $g_\mathscr{L}(\cdot \lvert \mu): \calM_1(E) \to \calM_1(E)$. Indeed, assume that $\nu \in \calM_1(E)$ is such that $g_\mathscr{L}(\mu | \nu)= 0$. Then, we necessarily have $\nu(x) > 0$ for all $x \in E$, and 
    by definition of $\Phi$ and the irreducibility of $\mathscr{L}$, the assumption that $g_{\mathscr{L}}(\nu | \mu) = 0$ implies that 
    \begin{align*}
        \frac{\mu(x)}{\nu(x)}\frac{\nu(y)}{\mu(y)} = 1 
    \end{align*}
    for all $x \neq y$. But this is equivalent to $\mu = \nu$.
    All in all, we have seen that 
    \begin{enumerate}[i.]
        \item $g_\mathscr{L}(\nu|\mu) \leq 0$ for all $\nu \in \calM_1(E)$ and 
        \item $g_\mathscr{L}(\nu | \mu) = 0$ if and only if $\nu = \mu$.
    \end{enumerate}
    When put together, these two properties imply that the functional
    \begin{align*}
        h(\cdot | \mu): \calM_1(E) \to \R, \quad \nu \mapsto h(\nu|\mu),
    \end{align*}
    is a strict Lyapunov function for the unique fixed point $\mu$ of the measure-valued ODE 
\begin{align}
    \partial_t \nu_t = \nu_t \mathscr{L}. 
\end{align}
Therefore, the fixed point $\mu$ is asymptotically stable and its basin of attraction is all of $\calM_1(E)$. But this is just another way of saying that $(\nu_t)_{t \geq 0}$ converges to $\mu$ as $t$ tends to infinity for all initial distributions $\nu \in \calM_1(E)$. 
    
This result is usually known as the \emph{ergodic theorem for finite-state Markov processes}. The proof given above shows that the convergence to the unique time-stationary measure also fits precisely into the physical picture of convergence to equilibrium. 
    
The rest of this paper is devoted to extending the results in this section to the setting of infinite-volume interacting particle systems. The philosophy of using the relative entropy as a Lyapunov functional will remain the same, but the proof itself becomes more technical. 

\section{Setting and results}\label{section:setting-and-results}
\subsection{Gibbs measures and interacting particle systems}

\subsubsection{Gibbs measures}
Let $q \in \N$ and consider $\Omega := \set{1,\dots, q}^{\Z^d}$. Equipped with the usual product topology and the corresponding Borel sigma-algebra $\calF$ this will serve as our configuration space. For $\Lambda \subset \Z^d$, let $\calF_{\Lambda}$ be the sub-sigma-algebra of $\calF$ that is generated by the open sets in $\Omega_{\Lambda}:=\set{1,\dots,q}^\Lambda$. We will use the shorthand notation $\Lambda \Subset \Z^d$ to signify that $\Lambda$ is a finite subset of $\Z^d$. In the following we will often denote for a given configuration $\omega \in \Omega$ by $\omega_{\Lambda}$ its projection to the volume $\Lambda \subset \Z^d$ and write $\omega_\Lambda\omega_\Delta$ for the finite-volume configuration in $\Lambda \cup \Delta$ composed of $\omega_\Lambda$ and $\omega_\Delta$ with disjoint $\Lambda, \Delta \Subset \Z^d$.  Denote the set of translation-invariant probability measures on $\Omega$ by $\calM_1^{inv}(\Omega)$. Then, for $\mu, \nu \in \calM_1^{inv}(\Omega)$ and a finite volume $\Lambda \Subset \Z^d$ define the relative entropy via 
\begin{align*}
    h_\Lambda(\nu | \mu) :=
    \begin{cases}
    \sum_{\omega_\Lambda \in \Omega_\Lambda}\nu(\omega_{\Lambda})\log\frac{\nu(\omega_{\Lambda})}{\mu(\omega_\Lambda)}, \quad &\text{if } \nu \ll \mu, 
    \\\
    \infty, &\text{else,}
    \end{cases}
\end{align*}
where we use the suggestive notation $\nu(\omega_{\Lambda}) = \nu(\set{\eta: \eta_{\Lambda} = \omega_\Lambda})$. Further, define the \textit{relative entropy density} via 
\begin{align*}
    h(\nu | \mu) := \limsup_{n \to \infty}\frac{1}{\abs{\Lambda_n}}h_{\Lambda_n}(\nu | \mu), 
\end{align*}
where $\Lambda_n := [-2^n +1, 2^n-1]^d$ is a sequence of hypercubes centered at the origin. 
Let us remark that in general one cannot replace the limes superior in the definition of the relative entropy density with a normal limit, even if both $\mu$ and $\nu$ are assumed to be translation invariant. However, if $\mu$ is a translation-invariant Gibbs measure, or even just asymptotically decoupled,  then one can show that the limit actually exists for all translation-invariant $\nu$. In that case, one can even replace the sequence of hypercubes by any increasing sequence $(\Delta_n)_{n \in \N}$ of subsets of $\Z^d$ that exhausts $\Z^d$ and satisfies the van Hove condition $\abs{\partial \Delta_n}/\abs{\Delta_n} \to 0$, where $\partial \Delta_n $ denotes the boundary of the set $\Delta_n$. We refer the interested reader to \cite[Chapter~2.6]{van_enter_regularity_1993} for more details.

We will be interested in situations where $\mu$ is a Gibbs measure for a translation-invariant non-null quasilocal specification on $\Omega$. 
\begin{definition}
A specification $\gamma = (\gamma_\Lambda)_{\Lambda \Subset \Z^d}$ is a family of probability kernels $\gamma_{\Lambda}$ from $\Omega _{\Lambda^c}$ to $\calM_1(\Omega)$ that additionally satisfies the following properties. 
\begin{enumerate}[i.]
    \item Each $\gamma_{\Lambda}$ is \emph{proper}, i.e., if $\Delta \subset \Lambda^c$, then 
    \begin{align*}
        \gamma_\Lambda(\eta_{\Lambda}\eta_\Delta | \eta_{\Lambda^c}) = \gamma_{\Lambda}(\eta_{\Lambda}|\eta_{\Lambda^c}) \mathbf{1}_{\eta_\Delta}(\eta_{\Lambda^c}). 
    \end{align*}
    \item The probability kernels are \emph{consistent} in the sense that if $\Delta \subset \Lambda \Subset \Z^d$, then 
    \begin{align*}
        \gamma_{\Lambda}(\gamma_{\Delta}(\eta_{\Delta}|\cdot)|\eta_{\Lambda^c}) 
        = 
        \gamma_\Lambda(\eta_\Delta | \eta_{\Lambda^c}). 
    \end{align*}
\end{enumerate}
\end{definition}

An infinite-volume probability measure $\mu$ on $\Omega$ is called a \textit{Gibbs measure} for $\gamma$ if $\mu$ satisfies the so-called \textit{DLR equations}, namely for all $\Lambda \Subset \Z^d$  and $\eta_\Lambda$ we have
\begin{align}\label{dlr-equations}
    \mu(\gamma_\Lambda(\eta_\Lambda | \cdot)) = \mu(\eta_\Lambda). 
\end{align}
We will denote the set of all Gibbs measures for a specification $\gamma$ by $\mathscr{G}(\gamma)$. 

For the existence and further properties of Gibbs measures one needs to impose some conditions on the specification $\gamma$. One sufficient condition for the existence of a Gibbs measure for a specification $\gamma$ is \textit{quasilocality}. 

\begin{definition}
A specification $\gamma$ is called 
\begin{enumerate}[i.]
    \item \emph{translation-invariant}, if for all $\Lambda \Subset \Z^d$ and $i \in \Z^d$ we have 
    \begin{align*}
        \gamma_{\Lambda +i}(\eta_{\Lambda +i}| \eta_{(\Lambda +i)^c}), 
    \end{align*}
    where $\Lambda+i$ denotes the lattice translate of $\Lambda$ by $i$.
    \item \emph{non-null}, if for some $\delta >0$
    \begin{align*}
        \inf_{\eta \in \Omega}\gamma_0(\eta_0 | \eta_{0^c}) \geq \delta.
    \end{align*}
    \item \emph{quasilocal}, if for all $\Lambda \Subset \Z^d$
    \begin{align*}
        \lim_{\Delta \uparrow \Z^d}\sup_{\eta, \xi \in \Omega}\abs{\gamma_{\Lambda}(\eta_\Lambda | \eta_{\Delta \setminus \Lambda}\xi_{\Delta^c}) - \gamma_\Lambda(\eta_{\Lambda}|\eta_{\Lambda^c})} = 0. 
    \end{align*}
\end{enumerate}
\end{definition}
We will sometimes consider the probability kernels $\gamma_\Lambda$ as functions $\Omega\to[0,1]$, $ \omega\mapsto\gamma_\Lambda(\omega_\Lambda | \omega_{\Lambda^c})$.
If $\gamma$ is a quasilocal specification, then each such map is then uniformly continuous.
For example, specifications defined via a translation-invariant uniformly absolutely summable potentials $\Phi=(\Phi_B)_{B \Subset \Z^d}$ are translation-invariant, non-null and quasilocal. Moreover, one can even show that, for Gibbs measures for such \textit{Gibbsian specifications}, the relative entropy density exists as a limit and not just as limes superior. For more details on Gibbs measures and specifications see~\cite{georgii_gibbs_2011},~\cite[Chapter 6]{friedli_statistical_2017} and~\cite[Chapter 4]{bovier_statistical_2006}. 

The DLR-formalism, which we used above, describes Gibbs measures through a collection of \textit{local} conditions~\eqref{dlr-equations}. For translation-invariant Gibbs measures there is also an alternative point of view that provides a \textit{global} description of Gibbs measures as the minimizers of a certain functional on the set $\calM_1^{inv}(\Omega)$. 
This is the content of the well-known \textit{Gibbs variational principle}. A detailed proof can be found in \cite[Chapter 6.9]{friedli_statistical_2017}. 
\begin{theorem}\label{theorem:gvp}
    Let $\Phi$ be an absolutely convergent and translation-invariant potential and  $\mu \in \mathscr{G}^{inv}(\Phi):=\mathscr{G}(\Phi)\cap \calM_1^{inv}(\Omega)$. 
    Define the relative entropy density by 
    \begin{align*}
        h(\nu |\Phi) 
        :=
        \limsup_{n \to \infty}\frac{1}{\abs{\Lambda_n}}h_{\Lambda_n}(\nu | \mu).
    \end{align*}
    Then,
    \begin{enumerate}[i.]
        \item For all $\nu \in \calM_1^{inv}(\Omega)$, $h(\nu |\Phi)$ exists as a limit and does not depend on $\mu$, only on $\Phi$,   
        \item $h(\nu |\Phi) \geq 0$ for all $\nu \in \calM_1^{inv}(\Omega)$ and 
        \item $h(\nu |\Phi) = 0$ if and only if $\nu \in \mathscr{G}^{inv}(\Phi)$. 
    \end{enumerate}
\end{theorem}

This \textit{static} global description is the starting point for our investigation. We aim for finding a \textit{dynamical} counterpart to the Gibbs variational principle, describing Gibbs measures as extremal points of another functional on $\calM_1^{inv}(\Omega)$, which describes the rate at which the relative entropy density changes when the system is subject to Markovian dynamics. 

\begin{remark}
    The Gibbs variational principle also provides an approach to define Gibbs measures for general measurable dynamical systems $(X, \mathcal{B}, T)$ that are additionally equipped with a potential $\phi:X \to \R$. 
    This is one possible starting point for the so-called \textit{thermodynamic formalism} for dynamical systems. 
    A non-trivial application of this thermodynamic formalism is the \textit{Bowen formula}, which can be used to calculate the Hausdorff dimension of attractors by constructing solutions to the analogue of the Gibbs variational problem in this situation. 
    An excellent reference is the monograph \cite{keller_equilibrium_1998}; other good resources are Ruelle's book \cite{ruelle_thermodynamic_2004} and the lecture notes \cite{bowen_equilibrium_1975}. An elementary proof of Bowen's formula in the simple situation of cookie-cutter maps can be found in \cite{falconer_techniques_1997}. 
\end{remark}

\subsubsection{Interacting particle systems}

We will consider time-continuous, translation-invariant Markov dynamics on $\Omega$, namely interacting particle systems characterized by time-homogeneous generators $\mathscr{L}$ with domain $\text{dom}(\mathscr{L})$ and its associated Markovian semigroup $(P_t)_{t \geq 0}$. 
For interacting particle systems we adopt the notation and exposition of the standard reference \cite[Chapter 1]{liggett_interacting_2005}. 
In our setting the generator $\mathscr{L}$ is given via a collection of translation-invariant transition rates $c_\Delta(\eta, \xi_\Delta)$, in finite volumes $\Delta \Subset \Z^d$, which are continuous in the starting configuration $\eta \in \Omega$. 
These rates can be interpreted as the infinitesimal rate at which the particles inside $\Delta$ switch from the configuration $\eta_\Delta$ to $\xi_\Delta$, given that the rest of the system is currently in state $\eta_{\Delta^c}$. 
The full dynamics of the interacting particle system is then given as the superposition of these local dynamics, i.e., 
\begin{align*}
    \mathscr{L}f(\eta) = \sum_{\Delta \Subset \Z^d}\sum_{\xi_\Delta}c_\Delta(\eta, \xi_\Delta)[f(\xi_\Delta \eta_{\Delta^c}) - f(\eta)].
\end{align*}
In \cite[Chapter 1]{liggett_interacting_2005} it is shown that the following two conditions are sufficient to guarantee the well-definedness. 
\begin{enumerate}[\bfseries (L1)]
    \item The total rate at which the particle at a particular site changes its spin is uniformly bounded, i.e.,
    \begin{align*}
        \sum_{\Delta \ni 0} \sum_{\xi_{\Delta}}\norm{c_{\Delta}(\cdot, \xi_{\Delta})}_{\infty} < \infty 
    \end{align*}
    \item and the total influence of a single coordinate on all other coordinates is uniformly bounded, i.e.,
    \begin{align*}
        \sum_{\Delta \ni 0}\sum_{x \neq 0}\sum_{\xi_{\Delta}}\delta_x\left(c_{\Delta}(\cdot, \xi_{\Delta})\right) < \infty, 
    \end{align*}
    where 
    \begin{align*}
        \delta_x(f) := \sup_{\eta, \xi: \ \eta_{x^c} = \xi_{x^c}}\abs{f(\eta)-f(\xi)}
    \end{align*}
    is the oscillation of a function $f:\Omega \to \R$ at the site $x$. 
\end{enumerate}
Under these conditions one can then show that the operator $\mathscr{L}$, defined as above, is the generator of a well-defined Markov process and that a core of $\mathscr{L}$ is given by 
\begin{align*}
    D(\Omega) := \Big\{ f \in C(\Omega): \quad \sum_{x \in \Z^d} \delta_x(f) < \infty\Big\}.
\end{align*}
Note that by considering the directional discrete derivatives  $\nabla^i_z : C(\Omega) \to C(\Omega)$ which are defined by
 \begin{align*}
     \nabla^i_z f(\eta) := f(\eta^{z,i}) - f(\eta), \quad f \in C(\Omega), i \in \set{1,\dots,q}, z \in \Z^d
 \end{align*}
with 
 \begin{align*}
     \eta^{z,i}_x
     :=
     \begin{cases}
      \eta_x, \quad &\text{if } x \neq z, \\\
      i, &\text{if } x = z, 
     \end{cases}
 \end{align*}
 then condition $\mathbf{(L2)}$ is equivalent to 
 \begin{align*}
     \sum_{\Delta \ni 0}\sum_{x \neq 0}\sum_{\xi_\Delta}\sum_{i=1}^q\norm{\nabla^i_x c_\Delta(\cdot, \xi_\Delta)}_\infty < \infty. 
 \end{align*}
  Let us note that the above conditions are given for the translation-invariant setting in which we exclusively work. They can be generalized to cover also non-translation-invariant models, for details see~\cite[Chapter 1]{liggett_interacting_2005}.
\subsection{Results}

Let us introduce some further conditions on the specification $\gamma = (\gamma_\Delta)_{\Delta \Subset \Z^d}$ and the rates $(c_\Delta(\cdot, \xi_\Delta))_{\Delta \Subset \Z^d, \xi_\Delta \in \Omega_\Delta}$ that will turn to be crucial for our results.

\medskip
\noindent
\textbf{Conditions for the specification.} 
\begin{enumerate}[\bfseries (S1)]
    \item $\gamma$ is quasilocal. 
    \item $\gamma$ is non-null.
    \item $\gamma$ satisfies 
    \begin{align*}
        \sum_{\Delta \ni 0:\  c_\Delta > 0}\sum_{z \neq 0}\sum_{i=1}^q \norm{\nabla^i_z \gamma_\Delta(\cdot | \cdot)}_\infty < \infty, 
    \end{align*}
    \item $\gamma$ is translation-invariant. 
\end{enumerate}
\noindent \textbf{Conditions for the rates.} 
\begin{enumerate}[\bfseries (R1)]
    \item For every $\Delta \Subset \Z^d$ and $\xi_\Delta \in \Omega_\Delta$ the function 
    \begin{align*}
        \Omega \ni \eta \mapsto c_\Delta(\eta, \xi_\Delta) \in [0,\infty)
    \end{align*}
    is continuous.  
    \item There are at most finitely many $\Delta \Subset \Z^d$ such that $0 \in \Delta$ and 
    \begin{align*}
        c_\Delta := \sup_{\xi_\Delta} \norm{c_\Delta(\cdot, \xi_\Delta)}_\infty > 0. 
    \end{align*}
    Denote by $R \in \N$ the maximal size of a subset $\Delta \Subset \Z^d$ with $c_\Delta >0$. 
    \item The total influence of all other particles on the particle at the origin is finite, i.e.,
    \begin{align*}
        \sum_{z \neq 0}\sum_{\Delta \ni 0}\sum_{\xi_\Delta}\sum_{i=1}^q\norm{\nabla^i_z c_\Delta(\cdot, \xi_\Delta)}_\infty < \infty. 
    \end{align*}
    \item The rates are translation-invariant, i.e., 
    \begin{align*}
        \forall x \in \Z^d \ \forall \Delta \Subset \Z^d \ \forall \eta \in \Omega: 
        \quad 
        c_{\Delta + x}(\tau_x \eta, \cdot) = c_\Delta(\eta, \cdot),
    \end{align*}
    where $\tau_x : \Omega \to \Omega$ is the lattice translation by $x$ acting on configurations. 
    \item The minimal transition rate is strictly positive, i.e.,
    \begin{align*}
        \inf_{\Delta \Subset \Z^d, \xi_{\Delta}, \eta: \ c_\Delta(\eta, \xi_\Delta)>0}c_\Delta(\eta, \xi_\Delta) > 0. 
    \end{align*}
\end{enumerate}
As previously noted, these assumptions guarantee that the interacting particle systems we consider are well-defined. 
For the last part of the dynamical Gibbs variational principle we will also need to assume irreducibility. 
\begin{enumerate}[\bfseries (R6)]
    \item The rates are irreducible, i.e., for every $\eta \in \Omega$, $\Delta \Subset \Z^d$ and $\xi_\Delta \in \Omega_\Delta$ there exists a finite sequence $\eta^{(0)}, \dots, \eta^{(n)}\subset\Omega$ such that $\eta^{(0)} = \eta$, $\eta^{(n)} = \xi_\Delta \eta_{\Delta^c}$ and the transition rate from $\eta^{(i)}$ to $\eta^{(i+1)}$ is positive for all $i \in {0, \dots, n-1}$. 
\end{enumerate}

\subsubsection{The (approximating) relative entropy loss}
Recall that $(P_t)_{t \geq 0}$ denotes the Markov semigroup corresponding to the Markov generator $\mathscr{L}$. We write $\nu_t:=\nu P_t$ for the time-evolved measure $\nu\in\calM_1(\Omega)$.
For $n \in \N$, the relative entropy loss in $\Lambda_n = [-2^n+1, 2^n -1]^d$ is defined by 
\begin{align*}
    g^n_{\mathscr{L}}(\nu|\mu) 
    := \frac{d}{dt}\lvert_{t=0}h_{\Lambda_n}( \nu_t\lvert \mu). 
\end{align*}
We define the \textit{relative entropy loss density} as
\begin{align*}
    g_{\mathscr{L}}(\nu\lvert \mu) := \limsup_{n \to \infty} \abs{\Lambda_n}^{-1}g^n_{\mathscr{L}}(\nu\lvert \mu). 
\end{align*}
In Proposition \ref{proposition:relative-entropy-loss-non-positive} we show that $g_{\mathscr{L}}(\nu\lvert \mu) \leq 0$, which justifies the name of $g_{\mathscr{L}}(\nu\lvert \mu)$. 
As it turns out, it will be more convenient to only consider the effect that spin flips inside the smaller cube 
$$\tilde{\Lambda}_n := [-2^n + n + 1, 2^n - n -1]^d
$$
have on the relative entropy. 
In Lemma~\ref{lemma:boundary-contributions-boundary-order} we will see that the contributions coming from spin flips at sites $x$ in the boundary region $\Lambda_n \setminus \tilde{\Lambda}_n$ are negligible in the density limit. 
This and the representation of $g^n_\mathscr{L}(\nu|\mu)$ in Lemma~\ref{lemma:finite-volume-relative-entropy-loss} motivate the definition of the \textit{approximating relative entropy loss in $\Lambda_n$} as
\begin{align*}
    \Tilde{g}^n_{\mathscr{L}}(\nu|\mu)
        := 
        \sum_{\eta_{\Lambda_n}}\sum_{\Delta \subset \Tilde{\Lambda}_n} \sum_{\xi_{\Delta}}
                \int \nu(d\omega)
                    c_{\Delta}(\omega, \xi_{\Delta})
                    \left[
                        \mathbf{1}_{\eta_{\Lambda_n}}(\xi_{\Delta}\omega_{\Delta^c}) 
                        -
                        \mathbf{1}_{\eta_{\Lambda_n}}(\omega)
                    \right]
                \log
            \left(
                \frac{\nu(\eta_{\Lambda_n})}{\mu(\eta_{\Lambda_n})}
            \right).
\end{align*}
The \textit{approximating relative entropy loss density} is then defined as 
\begin{align*}
    \Tilde{g}_\mathscr{L}(\nu | \mu) := \limsup_{n \to \infty}\abs{\Lambda_n}^{-1} \tilde{g}_\mathscr{L}^n(\nu | \mu). 
\end{align*}
\noindent 
Our first result relates the relative entropy loss density to the approximating relative entropy loss density and also shows that both quantities are non-positive. 

\begin{proposition}\label{proposition:relative-entropy-loss-non-positive}
    Suppose that the rates of an interacting particle system with generator $\mathscr{L}$ satisfy $\mathbf{(R1)}-\mathbf{(R5)}$. Moreover, assume that there exists a measure $\mu$ which is translation-invariant and time-stationary for the process generated by $\mathscr{L}$, such that $\mu \in \mathscr{G}(\gamma)$, where the specification $\gamma = (\gamma_{\Lambda})_{\Lambda \Subset \Z^d}$ satisfies $\mathbf{(S1)}-\mathbf{(S4)}$. 
    Then, for all $\nu \in \calM_1(\Omega)$, we have 
    \begin{align*}
            g_{\mathscr{L}}(\nu|\mu) 
            \leq 
            \tilde{g}_{\mathscr{L}}(\nu\lvert \mu)
            \leq 0. 
    \end{align*}
    In particular, we have for all $t\geq 0$,
    \begin{align*}
        h(\nu_t \lvert \mu) - h(\nu \lvert \mu) \leq 0. 
    \end{align*}
\end{proposition}

\noindent
The proof of Proposition \ref{proposition:relative-entropy-loss-non-positive} will be carried out in several steps and can be found after the proof of Lemma \ref{lemma:bulk-contribution}. 

\subsubsection{The dynamical Gibbs variational principle and the attractor property}

We are now ready to state our main result, generalizing the earlier works \cite{jahnel_attractor_2019}, \cite{kunsch_non_1984}, \cite{higuchi_results_1975}, and \cite{holley_free_1971}, to the setting of non-reversible interacting particle systems with finite local state spaces and updates in arbitrary finite regions. 

\begin{theorem}[Dynamical Gibbs variational principle]\label{theorem:dgvp}
    \sloppy Suppose that the rates $(c_{\Delta}(\cdot, \xi_{\Delta}))_{\Delta \Subset \Z^d, \xi_{\Delta}\in \Omega_{\Delta}}$ of an interacting particle system satisfy $\mathbf{(R1)}-\mathbf{(R6)}$.
    Moreover, assume that there exists a measure $\mu$ which is translation-invariant and time-stationary for the process generated by $\mathscr{L}$ such that $\mu \in \mathscr{G}(\gamma)$, where the specification $\gamma = (\gamma_{\Lambda})_{\Lambda \Subset \Z^d}$ satisfies $\mathbf{(S1)}-\mathbf{(S4)}$. 
    Then, the following variational principle for the approximating entropy loss density holds on the set of translation-invariant probability measures $\calM_1^{\text{inv}}(\Omega)$. 
    \begin{enumerate}[i.]
        \item For all $\nu \in \calM_1^{\text{inv}}(\Omega)$ the approximating entropy loss $\tilde{g}_{\mathscr{L}}(\nu \lvert \mu)$ exists as a limit in $\R \cup \set{-\infty}$. 
        \item The function $\nu \mapsto \tilde{g}_\mathscr{L}(\nu| \mu)$ is upper-semicontinuous on $\calM_1^{\text{inv}}(\Omega)$. 
        \item For all $\nu \in \calM_1^{\text{inv}}(\Omega)$ we have $\tilde{g}_\mathscr{L}(\nu \lvert \mu) \leq 0$. 
        \item For all $\nu \in \calM_1^{\text{inv}}(\Omega)$ with $\tilde{g}_\mathscr{L}(\nu \lvert \mu) = 0$ we have $\nu \in \mathscr{G}(\gamma)$. 
    \end{enumerate}
\end{theorem}
\noindent 

The proof of Theorem \ref{theorem:dgvp} will be carried out in several steps and can be found in Section \ref{section:proof-dgvp} after the proof of Lemma \ref{lemma:zero-loss-gibbs}. 
\begin{remark}
    Due to the irreducibility assumption $\mathbf{(R6)}$ on the rates, Theorem \ref{theorem:dgvp} does not  apply to the exclusion process, the voter model or the contact process. 
    While there is no hope to extend it to the voter model or the contact process, the exclusion process can be covered since it is irreducible when restricted to the subspaces
    \begin{align*}
        \Omega_\rho := \set{ \omega \in \Omega: \lim_{n\to \infty}\frac{\abs{\set{i \in \Lambda_n: \omega_i = 1}}}{\abs{\Lambda_n}} = \rho}, \quad \rho \in [0,1].  
    \end{align*}
    By being a bit more careful in the proof of Lemma \ref{lemma:zero-loss-gibbs} one can then indeed extend our result to $\mu, \nu$ as long as $\mu(\Omega_\rho) = 1 = \nu(\Omega_\rho)$ for some $\rho \in [0,1]$. 
\end{remark}

\begin{remark}
    One particular class of models to which our theory can be applied to are stochastic Ising models, if the specification $\gamma$ is defined via a translation-invariant potential $\Phi = (\Phi_B)_{B \Subset \Z^d}$ that satisfies 
    \begin{align*}
        \sum_{B \Subset \Z^d} \abs{B}\norm{\Phi_B}_{\infty} < \infty
    \end{align*}
    and the rates are of the general form 
    \begin{align*}
        c_\Delta(\eta, \xi_\Delta)
        =
        \begin{cases}
            \exp\left(-\beta\sum_{B\colon B \cap\Delta\neq\emptyset} \Phi_B(\xi_\Delta\eta_{\Delta^c})\right), \ \ &\text{if } \Delta = \set{x}, x\in \Z^d, \text{ and } \xi_x = - \eta_x. 
            \\\
            0, &\text{otherwise},
    \end{cases}
\end{align*}
    or even more general, with updates in larger regions with bounded diameter. 
    Then, the rates satisfy $\mathbf{(R1)}-\mathbf{(R6)}$ and the specification satisfies $\mathbf{(S1)}-\mathbf{(S4)}$ as one can see by using similar arguments as in the proof of \cite[Lemma 6.28]{friedli_statistical_2017}. 
    This class of models in particular includes such with long-range pairwise interactions $(\Phi_{\{x,y\}})_{x,y \in \Z^d}$ that satisfy
    \begin{align*}
        \norm{\Phi_{\{x,y\}}}_\infty \sim \norm{x-y}_1^{-\alpha}
    \end{align*} 
    for some $\alpha > d$.
\end{remark}

\noindent
 In the proof of Theorem \ref{theorem:dgvp}, we will see that on the set of non-null and translation-invariant measures we also have the same variational principle with the relative entropy loss density instead of the approximating relative entropy loss density. 
 
 \begin{corollary}\label{corollary:dgvp}
    Suppose that the rates $(c_{\Delta}(\cdot, \xi_{\Delta}))_{\Delta \Subset \Z^d, \xi_{\Delta}\in \Omega_{\Delta}}$ of an interacting particle system satisfy $\mathbf{(R1)}-\mathbf{(R6)}$.
    Moreover, assume that there exists a measure $\mu$ which is translation-invariant and time-stationary for the process generated by $\mathscr{L}$, such that $\mu \in \mathscr{G}(\gamma)$, where the specification $\gamma$ satisfies $\mathbf{(S1)}-\mathbf{(S4)}$.
    Then, the following variational principle for the relative entropy loss density holds on the set of non-null translation-invariant probability measures $\calM_1^{\text{inv,non-null}}(\Omega)$. 
    \begin{enumerate}[i.]
        \item For all $\nu \in \calM_1^{\text{inv,non-null}}(\Omega)$ the approximating entropy loss $g_\mathscr{L}(\nu \lvert \mu)$ exists as a limit in $\R \cup \set{-\infty}$. 
        \item The function $\nu \mapsto g_\mathscr{L}(\nu|\mu)$ is upper-semicontinuous on $\calM_1^{\text{inv,non-null}}(\Omega)$. 
        \item For all $\nu \in \calM_1^{\text{inv,non-null}}(\Omega)$ we have $g_\mathscr{L}(\nu \lvert \mu) \leq 0$. 
        \item For all $\nu \in \calM_1^{\text{inv.non-null}}(\Omega)$ with $g_\mathscr{L}(\nu \lvert \mu) = 0$ we have $\nu \in \mathscr{G}(\gamma)$. 
    \end{enumerate}
 \end{corollary}
 
The main difference with the results in \cite{jahnel_attractor_2019} is conclusion $iv.$ in Theorem~\ref{theorem:dgvp} and Corollary~\ref{corollary:dgvp}. For previous versions of attractor properties for irreversible systems \cite[Theorem 2.8 and Theorem 2.10]{jahnel_attractor_2019}, this property of the (approximating) relative entropy loss had to be added as an extra assumption. Only in the reversible setting the techniques in \cite{jahnel_attractor_2019} were able to show that this property holds in fairly general situations. 
By using finer estimates and an extension of a powerful connection between the rates of the forward process and a suitable backward process, which was first used in a similar context in \cite{kunsch_non_1984} and will be discussed in more detail in Section \ref{section:proof-time-reversal}, we were now able to replace the reversibility assumption by the much weaker assumption of time-stationarity.

We started out with the big goal of investigating the set of possible limit points of $(\nu_t)_{t \geq 0}$ for translation-invariant initial measures $\nu_0$ and non-reversible interacting particle systems. 
With the dynamical Gibbs variational principle at hand, it is now easy to show that all possible limit points are themselves translation-invariant Gibbs measures with respect to the same specification. This extends the attractor properties in \cite{jahnel_attractor_2019} to non-reversible systems, without adding further assumptions like \cite[Condition 2.5 and Condition 2.9]{jahnel_attractor_2019}.

\begin{theorem}[Attractor property for irreversible interacting particle systems]\label{theorem:attractor-property}
    \sloppy Assume that the rates $(c_{\Delta}(\cdot, \xi_{\Delta}))_{\Delta \Subset \Z^d,\xi_\Delta\in \Omega_\Delta}$ satisfy conditions $\mathbf{(R1)}-\mathbf{(R6)}$. Moreover, assume that there exists a translation-invariant time-stationary measure $\mu$ which is a Gibbs measure with respect to a specification $\gamma$ that satisfies $\mathbf{(S1)}-\mathbf{(S4)}$. 
    Then, the $\omega$-limit set of the family of translation-invariant probability measures $\calM_1^{inv}(\Omega)$ is $\mathscr{G}(\gamma)$, 
    i.e., for any translation-invariant starting measure $\nu\in \calM_1^{inv}(\Omega)$ where the sequence $(\nu_{t_n})_{n \in \N}$ 
    converges weakly to $\nu^*$ as $t_n \uparrow \infty$, we have that $\nu^* \in \mathscr{G}(\gamma)$. 
\end{theorem}

The proof of Theorem ~\ref{theorem:attractor-property} can be found at the end of Section ~\ref{section:proof-attractor-property}. 
\subsubsection{An alternative characterization of time-stationary measures}

In the proof of Theorem \ref{theorem:dgvp} we will encounter an auxiliary process with rates given by
\begin{align*}
    \hat{c}_{\Delta}(\eta, \xi_{\Delta}) := c_{\Delta}(\xi_{\Delta}\eta_{\Delta^c}, \eta_{\Delta})\frac{\gamma_{\Delta}(\xi_{\Delta}|\eta_{\Delta^c})}{\gamma_{\Delta}(\eta_{\Delta}|\eta_{\Delta^c})},
\end{align*}
which can be interpreted as the time-reversal of the original process w.r.t.~to the stationary measure $\mu$ with local conditional distributions given by $\gamma$. 
In the reversible case, the rates of the time-reversed process agree with the original rates, which is known as the detailed balance equation. In the irreversible case, this does not hold, but we will see that time-stationarity implies that the weaker condition 
\begin{align}\label{oscillation-eqn}
    \sum_{\Delta \Subset \Z^d}\sum_{\xi_{\Delta}}\nabla_{\Lambda}\left(c_{\Delta}(\cdot, \xi_{\Delta})- \hat{c}(\cdot, \xi_{\Delta})\right) (\eta) = 0 
\end{align}
holds for all $\Lambda\Subset\Z^d$ and $\eta \in \Omega$, see Proposition \ref{proposition:oscillation-equation}, where $\nabla_\Lambda$ is the generalized differential operator defined by 
\begin{align*}
    \nabla_\Lambda f(\eta) := \sum_{\xi_\Lambda}[f(\xi_\Lambda \eta_{\Lambda^c})-f(\eta)]. 
\end{align*}
It is well known that the detailed balance equations are equivalent to the reversibility of the measure, and it is natural to ask whether equation \eqref{oscillation-eqn} is equivalent to the time-stationarity of a measure $\mu$ with local conditional distributions given by $\gamma$. We show that this is true under an additional assumption on the mixing coefficients of the measure $\mu$. 
For a measure $\nu$ and a subvolume $\Lambda \Subset \Z^d$ we define the \textit{mixing coefficients}
\begin{align*}
    \alpha_{\nu}(\Lambda, n) := \sup\set{\abs{\nu(A\cap B)-\nu(A)\nu(B)}\colon A \in \calF_{\Lambda}, B \in \calF_{\Delta_n^c}}, \quad n \in \N, 
\end{align*}
where $\Delta_n := [-n,n]^d$. Roughly speaking, this measures how much the spins inside the finite volume $\Lambda$ are correlated with the spins outside the box $\Delta_n$. 
We additionally need to assume that the  convergence of the sums in \eqref{oscillation-eqn} is uniform in $z \in \Z^d$.
For this we define for $n \in \N$
\begin{align*}
    \beta(n) 
    := 
    \sup_{z \in \Z^d, i \in \set{1,\dots,q}} 
    \max \left( 
        \sum_{\Delta \cap \Delta_n = \emptyset}\sum_{\xi_\Delta}\norm{\nabla_z^i c_\Delta(\cdot, \xi_\Delta)}_\infty, 
        \sum_{\Delta \cap \Delta_n = \emptyset}\norm{\nabla_z^i \hat{c}_\Delta(\cdot, \xi_\Delta)}_\infty
    \right). 
\end{align*}

\begin{theorem}[A mixing criterion for time-stationarity]\label{theorem:all-gibbs-invariant}
    \sloppy
    Assume that the rates $(c_\Delta(\cdot, \xi_\Delta))_{\Delta \Subset \Z^d, \xi_\Delta\in \Omega_\Delta}$ satisfy conditions $\mathbf{(R1)}-\mathbf{(R4)}$ and the specification $\gamma$ satisfies the conditions $\mathbf{(S1)}-\mathbf{(S4)}$,
    Moreover, assume that we have 
    \begin{align*}
        \beta(0) < \infty \quad \text{and} \quad \beta(n) \to 0 \ \text{as} \ n  \to \infty,
    \end{align*}
    and that for all $\Lambda \Subset \Z^d$ and $\eta \in \Omega$ we have 
    \begin{align*}
         \sum_{\Delta \Subset \Z^d}\sum_{\xi_{\Delta}}\nabla_{\Lambda}\left(c_{\Delta}(\cdot, \xi_{\Delta})- \hat{c}(\cdot, \xi_{\Delta})\right) (\eta) = 0. 
    \end{align*}
    If the mixing coefficients of $\mu \in \mathscr{G}(\gamma)$ satisfy 
    \begin{align*}
        \forall \Lambda \Subset \Z^d: \quad \sum_{n \in \N}\alpha_{\mu}(\Lambda,n)n^{d-1} < \infty,
    \end{align*}
    then $\mu$ is time-stationary for the interacting particle system associated to the rates $(c_\Delta(\cdot, \xi_\Delta))_{\Delta \Subset \Z^d, \xi_\Delta\in \Omega_\Delta}$. 
\end{theorem}

The proof of this result can be found after the proof of Proposition~\ref{proposition:oscillation-equation} at the end of Section~\ref{section:proof-time-reversal}. 

Let us note that the mixing condition is hard to verify in practice, especially in non-uniqueness regimes. 
An example of a situation where the mixing condition can be verified is Dobrushin's uniqueness regime, see~\cite[Chapter 8]{georgii_gibbs_2011} and the discussion in \cite{kunsch_non_1984}. We can imagine that the mixing condition can also be verified in the regime of Gaussian concentration~\cite{chazottes_gaussian_2020,chazottes_relative_2022}. 

\section{Proofs}\label{section:proofs}
\subsection{Proof strategy}
The proof of Theorem \ref{theorem:dgvp} proceeds in several steps. We start by deriving an explicit formula for the finite-volume relative entropy loss in terms of the generator and the time-stationary measure $\mu$. We will then show that the relative entropy density is non-increasing along trajectories by rewriting it similarly as in~\eqref{eqn:relative-entropy-loss-finite-state-spaces}. 
However, this is not as straightforward as in the finite-volume case and we need to find appropriate replacements for our finite-volume arguments.
For motivational purposes, let us briefly go back to finite state spaces. Consider an interacting particle system with irreducible transition rates $(c_\Delta(\cdot, \xi_\Delta))_{\Delta \subset \Lambda, \xi_\Delta \in \Omega_\Delta}$ in a finite volume $\Lambda \Subset \Z^d$.
In Section~\ref{section:finite-state-space}, we used that a probability measure $\mu$ is time-stationary w.r.t.~the dynamics if and only if it satisfies~\eqref{proof:sign-relative-entropy-loss-finite-space-1}. 
In the situation we consider here, this equation takes the form
\begin{align*}
    \forall \eta_\Lambda \in \Omega_\Lambda: 
    \quad 
    \sum_{\Delta \subset \Lambda}\sum_{\xi_\Delta}\mu(\eta_\Lambda) c_\Delta(\eta_\Lambda, \xi_\Delta) 
    =
    \sum_{\Delta \subset \Lambda}\sum_{\xi_\Delta}\mu(\xi_\Delta \eta_{\Lambda \setminus \Delta})c_\Delta(\xi_\Delta\eta_{\Lambda \setminus \Delta}, \eta_\Delta),
\end{align*}
which we can rewrite as 
\begin{align}\label{time-reversal-motivation-1}
    \forall \eta_\Lambda \in \Omega_\Lambda: 
    \quad 
    \sum_{\Delta \subset \Lambda}\sum_{\xi_\Delta}\left(c_\Delta(\eta_\Lambda, \xi_\Delta) - \frac{\mu(\xi_\Delta \eta_{\Lambda \setminus \Delta})}{\mu(\eta_\Lambda)}c_\Delta(\xi_\Delta\eta_{\Lambda \setminus \Delta}, \eta_\Delta)\right)
    =
    0. 
\end{align}
In a way to be made precise, the terms 
$$
\hat{c}_\Delta(\eta, \xi_\Delta) 
:=
\frac{\mu(\xi_\Delta \eta_{\Lambda \setminus \Delta})}{\mu(\eta_\Lambda)}c_\Delta(\xi_\Delta\eta_{\Lambda \setminus \Delta}, \eta_\Delta)
$$ 
are the rates of the \textit{time-reversed} version of the interacting particle system we consider. 
However, for general interacting particle systems, in infinite volumes the equation \eqref{time-reversal-motivation-1} does not make sense and we can not hope to use it naively in our quest to derive an analogue of  \eqref{eqn:relative-entropy-loss-finite-state-spaces} for irreversible interacting particle systems in infinite volumes. 
As it turns out, in finite volumes it already suffices to know that 
\begin{align}\label{motivation-oscillation-equation}
    \forall \eta_\Lambda \in \Omega_\Lambda \ \forall z \in \Lambda \ \forall i \in  \set{1,\dots,q}: 
    \quad 
    \sum_{\Delta \subset \Lambda}\sum_{\xi_\Delta}\nabla^i_z\left(c_\Delta(\cdot, \xi_\Delta) - \hat{c}_\Delta(\cdot, \xi_\Delta) \right)(\eta_\Lambda) = 0.
\end{align}

The idea of using~\eqref{motivation-oscillation-equation} to extend the works of Holley, and Higuchi and Shiga to non-reversible systems was first used in \cite{kunsch_non_1984}. Section \ref{section:proof-time-reversal} is devoted to extending this work to our much more general setting. Whereas Künsch only considered interacting particle systems on a binary local state space with single-site updates, we allow general finite local state spaces and updates in arbitrary finite regions. For this reason, we cannot rely on any spin-flip symmetry arguments and need to find appropriate replacements for the corresponding steps.  We will therefore first establish an infinite-volume analogue of~\eqref{motivation-oscillation-equation} and then use it to prove Proposition~\ref{proposition:relative-entropy-loss-non-positive} in Section~\ref{section:proof-entropy-loss-non-positive}. 

To establish the remaining parts of Theorem \ref{theorem:dgvp}, i.e., 
that the approximating relative entropy loss exists as a limit and that this limit defines an upper-semicontinuous functional on the space of translation-invariant probability measures, is then our next main task. We again have to eliminate dangerous terms and use a subadditivity argument to get the claimed convergence and semicontinuity. Equipped with these intermediate results, we are then ready to show the last step in the dynamical Gibbs variational principle, namely that translation-invariant measures $\nu$ with $\Tilde{g}_\mathscr{L}(\nu|\mu)=0$ are also Gibbs measures compatible with $\gamma$.
By putting all of our previous results together, we are then ready to prove the attractor property for non-reversible interacting particle systems in infinite volumes. 

\subsection{The time-reversal rates and the oscillation equations}\label{section:proof-time-reversal}

We start with an elementary integral identity on which we will rely heavily for the rest of this subsection. 

\begin{lemma}[Switching Lemma]\label{generalized-switching-lemma}
    Let $\gamma = (\gamma_{\Lambda})_{\Lambda \Subset \Z^d}$ be a specification, $\mu \in \mathscr{G}(\gamma)$ and $(c_{\Delta}(\cdot, \xi_\Delta))_{\Delta \Subset \Z^d, \xi_\Delta \in \Omega_\Delta}$ the rates of an interacting particle system. 
    Additionally, assume that $\gamma$ is strictly positive, i.e., that we have 
    \begin{align*}
    \gamma_{\Lambda}(\eta_\Lambda|\eta_{\Lambda^c}) > 0
    \end{align*}
    for all $\Lambda \Subset \Z^d$ and $\eta \in \Omega$. 
    Then, for all bounded and measurable $f,g : \Omega \to \R$ and $\Delta \Subset \Z^d$ we have 
    \begin{align}\label{generalized-switching-lemma-identity}
        \sum_{\xi_{\Delta} \in \Omega_{\Delta}}\int_{\Omega}c_{\Delta}(\omega, \xi_{\Delta})f(\omega)g(\xi_{\Delta}\omega_{\Delta^c})\mu(d\omega) 
        = 
        \sum_{\xi_{\Delta}\in \Omega_{\Delta}}\int_{\Omega}\hat{c}_{\Delta}(\omega, \xi_{\Delta})f(\xi_{\Delta}\omega_{\Delta^c})g(\omega)\mu(d\omega), 
    \end{align}
    where
    \begin{align}
        \hat{c}_{\Delta}(\eta, \xi_{\Delta}) := c_{\Delta}(\xi_{\Delta}\eta_{\Delta^c}, \eta_{\Delta})\frac{\gamma_{\Delta}(\xi_{\Delta}|\eta_{\Delta^c})}{\gamma_{\Delta}(\eta_{\Delta}|\eta_{\Delta^c})}. 
    \end{align}
\end{lemma}

\noindent 

For simplicity, in the following, we will sometimes denote integration w.r.t.~$\mu$ by $\bbE[\cdot]$. 

\begin{proof}
    As a first step, note that, for fixed $\Delta \Subset \Z^d$ and $\xi_\Delta \in \Omega_\Delta$, the maps 
    \begin{align*}
        \Omega \ni \omega \mapsto g(\xi_{\Delta}\omega_{\Delta^c}) \in \R, \quad \Omega \ni \omega \mapsto f(\xi_{\Delta}\omega_{\Delta^c}) \in \R,
    \end{align*}
    are $\calF_{\Delta^c}$-measurable. 
    Therefore, we can use that $\gamma$ is a version of the local conditional distribution of $\mu$ and the definition of the rates $\hat{c}$ to obtain the $\mu$-almost sure identity 
    \begin{align*}
        \condexp{c_{\Delta}(\cdot, \xi_{\Delta})f(\cdot)g(\xi_{\Delta}\cdot_{\Delta^c})}{\calF_{\Delta^c}}(\omega) 
        &=
        g(\xi_{\Delta}\omega_{\Delta^c}) \condexp{c_{\Delta}(\cdot, \xi_{\Delta})f(\cdot)}{\calF_{\Delta^c}}(\omega)
        \nonumber \\\
        &=
        g(\xi_{\Delta}\omega_{\Delta^c}) 
        \sum_{\zeta_{\Delta}}
            \gamma_{\Delta}(\zeta_{\Delta}|\omega_{\Delta^c})c_{\Delta}(\zeta_{\Delta}\omega_{\Delta^c}, \xi_{\Delta})f(\zeta_{\Delta}\omega_{\Delta^c})
        \\\
        &= 
        g(\xi_{\Delta}\omega_{\Delta^c}) 
        \sum_{\zeta_{\Delta}}
            \gamma_{\Delta}(\xi_{\Delta}|\omega_{\Delta^c})\hat{c}_{\Delta}(\xi_{\Delta}\omega_{\Delta^c}, \xi_{\Delta})f(\zeta_{\Delta}\omega_{\Delta^c}). 
        \nonumber 
    \end{align*}
    If we now sum this over $\xi_{\Delta} \in \Omega_{\Delta}$, exchange the order of summation and apply the same arguments as above in reverse – with $f$ taking the role of $g$ and vice versa – we get
    \begin{align*}
        \sum_{\xi_{\Delta}}
            \condexp{c_{\Delta}(\cdot, \xi_{\Delta})f(\cdot)g(\xi_{\Delta}\cdot_{\Delta^c})}{\calF_{\Delta^c}}(\eta) 
        = 
        \sum_{\zeta_{\Delta}}
            \condexp{\hat{c}_{\Delta}(\cdot, \zeta_{\Delta})f(\zeta_{\Delta}\cdot_{\Delta^c})g(\cdot)}{\calF_{\Delta^c}}(\eta).
    \end{align*}
    By integrating both sides with respect to $\mu$ and applying the law of total expectation, we obtain 
    \begin{align*}
        \sum_{\xi_{\Delta}}\int_{\Omega}c_{\Delta}(\omega, \xi_{\Delta})f(\omega)g(\xi_{\Delta}\omega_{\Delta^c})\mu(d\omega) 
        =
        \sum_{\zeta_{\Delta}}\int_{\Omega}\hat{c}_{\Delta}(\omega, \zeta_{\Delta})f(\zeta_{\Delta}\omega_{\Delta^c})g(\omega)\mu(d\omega),
    \end{align*}
    which completes the proof.
\end{proof}

We will often have to estimate terms where the specification appears in the denominator. The main tool for obtaining bounds will be the following lemma. 

\begin{lemma}\label{lemma:non-nullness-inequality}
    Let $\Delta \subset \Lambda \Subset \Z^d$ and $\rho$ be a probability measure that is non-null with parameter $\delta(\rho) >0$, then for all $\eta, \xi \in \Omega$ we have
    \begin{align}\label{log-bound-non-nullness}
        \abs{\log \left( \frac{\rho(\xi_{\Delta}\eta_{\Lambda \setminus \Delta})}{\rho(\eta_{\Lambda})}\right)}
        =
        \abs{\log \left( \frac{\rho(\xi_{\Delta}\lvert \eta_{\Lambda \setminus \Delta})}{\rho(\eta_{\Delta}\lvert \eta_{\Lambda \setminus \Delta})}\right)} 
        \leq 
        \abs{\Delta}\log \left( \frac{1}{\delta(\rho)}\right).
    \end{align}
    In particular, for $\rho \in \mathscr{G}(\gamma)$ with non-null specification $\gamma$, the same estimate holds. 
\end{lemma}
This estimate already appears in \cite{jahnel_attractor_2019}, but for the sake of being self-contained we also give the short proof here. 
\begin{proof}
    The first identity is clear by definition of conditional probabilities, so we only have to show the inequality. For this, fix an enumeration $i_1, \dots, i_k$ of the elements of $\Delta$ and introduce the notation 
    \begin{align*}
        [i_j, i_k] := \set{i_j, i_{j+1},\dots, i_k}, \quad 1 \leq j \leq k. 
    \end{align*}
    With this at hand, we can use the chain rule for conditional probabilities to write
    \begin{align}\label{proof:chain-rule}
        \rho(\xi_{\Delta}\lvert \eta_{\Lambda \setminus \Delta})
        =
        \prod_{j=1}^{k-1}
            \rho(\eta_{i_j}| \eta_{[i_{j+1},i_k]}\eta_{\Lambda \setminus \Delta}).
    \end{align}
    Now we show that each factor is bounded from below by $\delta$. Indeed, via an elementary calculation, we see that
    \begin{align*}
        \rho(\eta_{i_j}| \eta_{[i_{j+1},i_k]}&\eta_{\Lambda \setminus \Delta})
        =
        \frac{
            \int\rho(d\omega)\rho(\eta_{[i_j, i_k]}\eta_{\Lambda \setminus \Delta}| \omega_{\Lambda^c \cup [i_1, i_{j-1}]})
        }
        {
            \int\rho(d\omega)\rho(\eta_{[i_{j+1}, i_k]}\eta_{\Lambda \setminus \Delta}| \omega_{\Lambda^c \cup [i_1, i_{j-1}]})
        }
        \\\
        &= 
        \frac{
            \int
            \frac{\rho(d\omega)\rho(\eta_{[i_j, i_k]}\eta_{\Lambda \setminus \Delta}| \omega_{\Lambda^c \cup [i_1, i_{j-1}]})}
            {\rho(\eta_{[i_{j+1}, i_k]}\eta_{\Lambda \setminus \Delta}| \omega_{\Lambda^c \cup [i_1, i_{j-1}]})}
            \rho(\eta_{[i_{j+1}, i_k]}\eta_{\Lambda \setminus \Delta}| \omega_{\Lambda^c \cup [i_1, i_{j-1}]})
        }
        {
            \int\rho(d\omega)\rho(\eta_{[i_{j+1}, i_k]}\eta_{\Lambda \setminus \Delta}| \omega_{\Lambda^c \cup [i_1, i_{j-1}]})
        }
        \\\
        &= 
        \frac{
            \int \rho(d\omega) \rho(\eta_{i_j}| \eta_{[i_{j+1}, i_k]}\eta_{\Lambda \setminus \Delta}\omega_{\Lambda^c\cup [i_1, i_{j-1}]})
            \rho(\eta_{[i_{j+1}, i_k]}\eta_{\Lambda \setminus \Delta}| \omega_{\Lambda^c \cup [i_1, i_{j-1}]})
        }
        {
            \int\rho(d\omega)\rho(\eta_{[i_{j+1}, i_k]}\eta_{\Lambda \setminus \Delta}| \omega_{\Lambda^c \cup [i_1, i_{j-1}]})
            }\geq 
        \delta.
    \end{align*}
    In conjunction with the representation~\eqref{proof:chain-rule}, this implies the desired upper bound. 
    If $\rho \in \mathscr{G}(\gamma)$ for a non-null specification $\gamma$, then we can carry out exactly the same calculations as before, except that we need to use the DLR equations to write
    \begin{align*}
        \rho(\eta_{i_j}| \eta_{[i_{j+1}, i_k]}\eta_{\Lambda \setminus \Delta}\omega_{\Lambda^c\cup [i_1, i_{j-1}]})
        =
        \gamma_{i_j}(\eta_{i_j}| \eta_{[i_{j+1}, i_k]}\eta_{\Lambda \setminus \Delta}\omega_{\Lambda^c\cup [i_1, i_{j-1}]}). 
    \end{align*}
    This finishes the proof.
\end{proof}

As a first step, we now verify that the regularity of the original rates and the specification also implies some regularity for the time-reversal rates.  This technical calculation will not only be needed to ensure that the infinite sums occurring in Proposition \ref{proposition:oscillation-equation} are well-defined, but we will also use it to show that the relative entropy is non-increasing in the proof of  Theorem \ref{theorem:dgvp}.

\begin{lemma}\label{lemma:regularity-time-reversal}
    Assume that the rates $(c_{\Delta}(\cdot, \xi_{\Delta}))_{\Delta \Subset \Z^d, \xi_{\Delta}\in \Omega_{\Delta}}$ of an interacting particle system satisfy the conditions 
    $\mathbf{(R1)}-\mathbf{(R4)}$ and that there exists a measure $\mu$ which is time-stationary for the process generated by $\mathscr{L}$ and such that $\mu \in \mathscr{G}(\gamma)$, where the specification $\gamma = (\gamma_{\Delta})_{\Delta \Subset \Z^d}$ satisfies the conditions $\mathbf{(S1)}-\mathbf{(S4)}$. 
    Then, the time-reversal rates $(\hat{c}_\Delta(\cdot, \xi_\Delta))_{\Delta \Subset \Z^d, \xi_\Delta \in \Omega_\Delta}$ have the following properties.  
    \begin{enumerate}[i.]
        \item The total rate of change of a single site is uniformly bounded over all sites, i.e.,
        \begin{align*}
            \sup_{x \in \Z^d} \sum_{\Delta \ni x}\sum_{\xi_{\Delta}}\norm{\hat{c}_\Delta(\cdot, \xi_\Delta)}_\infty < \infty. 
        \end{align*}
        \item The total influence of all other sites on a fixed site is uniformly bounded over all sites, i.e.,
        \begin{align*}
            \sup_{y \in \Z^d}\sum_{\Delta \ni y}\sum_{z \neq y}\sum_{\xi_{\Delta}}\sum_{i=1}^q \norm{\nabla^i_z \hat{c}_\Delta(\cdot, \xi_\Delta)}_\infty < \infty. 
        \end{align*}
        \item For all $z \in \Z^d$ it holds that
        \begin{align*}
            \sum_{\Delta \Subset \Z^d}\sum_{\xi_\Delta}\sum_{i=1}^q \norm{\nabla^i_z \hat{c}_\Delta(\cdot , \xi_\Delta)}_\infty < \infty. 
        \end{align*}
        \item The time-reversal rates $(\hat{c}_\Delta(\cdot, \xi_\Delta))_{\Delta \Subset \Z^d, \xi_\Delta \in \Omega_\Delta}$ are also translation-invariant.
    \end{enumerate}
\end{lemma}
Before we give the proof, note that the translation invariance of the rates and the specification implies that $\mathbf{(R3)}$ and $\mathbf{(S3)}$ are respectively equivalent to 
\begin{align*}
            \forall z \in \Z^d: \quad \sum_{\Delta \Subset \Z^d}\sum_{\xi_\Delta}\sum_{i=1}^q \norm{\nabla^i_z c_\Delta(\cdot , \xi_\Delta)}_\infty < \infty,
\end{align*}
and 
\begin{align*}
    \forall z \in \Z^d: \quad \sum_{\Delta \Subset \Z^d}\sum_{i=1}^q \norm{\nabla^i_z \gamma_\Delta(\cdot | \cdot)}_\infty < \infty, 
\end{align*}
as one can easily verify by elementary calculations. 
\begin{proof}
\textit{Ad i.: } 
For fixed $\Delta \Subset \Z^d$, $\xi_\Delta \in \Omega_\Delta$ and $\eta \in \Omega$ we have by Lemma \ref{lemma:non-nullness-inequality} and assumption $(\mathbf{R2})$
\begin{align*}
    \abs{\hat{c}_\Delta(\eta, \xi_\Delta)} 
    =
    \abs{c_\Delta(\xi_\Delta \eta_{\Delta^c}, \eta_\Delta) \frac{\gamma_\Delta(\xi_{\Delta}| \eta_{\Delta^c})}{\gamma_\Delta(\eta_\Delta| \eta_{\Delta^c})}}
    \leq 
    \frac{1}{\delta}e^{\abs{\Delta}} \norm{c_\Delta(\cdot, \eta_{\Delta})}_\infty \leq \frac{1}{\delta}e^R \norm{c_\Delta(\cdot, \eta_\Delta)}_\infty. 
\end{align*}
By assumptions $\mathbf{(R1)},\mathbf{(R2)}$, and $\mathbf{(R4)}$ we have 
\begin{align*}
    \sup_{x \in \Z^d} \sum_{\Delta \ni x} \sum_{\xi_{\Delta}xe}\norm{c_{\Delta}(\cdot, \xi_{\Delta})}_{\infty} < \infty, 
\end{align*}
and this implies 
\begin{align*}
    \sup_{x \in \Z^d} \sum_{\Delta \ni x}\sum_{\xi_\Delta}\norm{\hat{c}_\Delta(\cdot, \xi_{\Delta})}_\infty 
    \leq 
    \sup_{x \in \Z^d}\sum_{\Delta \ni x} \sum_{\eta_\Delta}\frac{1}{\delta}e^R \norm{c_\Delta(\cdot, \eta_\Delta)}_\infty < \infty. 
\end{align*}
\textit{Ad ii.: }
For fixed $z \in \Z^d$ and $i \in \set{1,\dots, q}$ we have 
\begin{align*}
        \abs{\hat{c}_{\Delta}(\eta^{z,i}, \xi_{\Delta})
        -
        \hat{c}_{\Delta}(\eta, \xi_{\Delta})}
        =
        &\abs{
            c_{\Delta}(\xi_{\Delta}\eta^{z,i}_{\Delta^c}, \eta^{z,i}_{\Delta}) 
            \frac{\gamma_{\Delta}(\xi_{\Delta}|\eta^{z,i}_{\Delta^c}) }{\gamma_{\Delta}(\eta^{z,i}_{\Delta}|\eta^{z,i}_{\Delta^c})}
            -
            c_{\Delta}(\xi_{\Delta}\eta_{\Delta^c}, \eta_{\Delta}) 
            \frac{\gamma_{\Delta}(\xi_{\Delta}|\eta_{\Delta^c}) }{\gamma_{\Delta}(\eta_{\Delta}|\eta_{\Delta^c})}
        }
        \\\
        \leq 
        &\abs{c_{\Delta}(\xi_{\Delta}\eta^{z,i}_{\Delta^c}, \eta_{\Delta}^{z,i})}
        \abs{
            \frac{\gamma_{\Delta}(\xi_{\Delta}|\eta^{z,i}_{\Delta^c}) }{\gamma_{\Delta}(\eta^{z,i}_{\Delta}|\eta^{z,i}_{\Delta^c})}
            -
            \frac{\gamma_{\Delta}(\xi_{\Delta}|\eta_{\Delta^c}) }{\gamma_{\Delta}(\eta_{\Delta}|\eta_{\Delta^c})}
        }
        \\\
        +
        &\abs{\frac{\gamma_{\Delta}(\xi_{\Delta}|\eta_{\Delta^c}) }{\gamma_{\Delta}(\eta_{\Delta}|\eta_{\Delta^c})}}
        \abs{c_{\Delta}(\xi_{\Delta}\eta^{z,i}_{\Delta^c}, \eta_{\Delta}^{z,i})-c_{\Delta}(\xi_{\Delta}\eta_{\Delta^c}, \eta_{\Delta})}.
\end{align*}
To estimate this further, we will have to make a case distinction over whether the site $z$ is contained in $\Delta$ or not. If $z$ is contained in $\Delta$, then we can naively use Lemma \ref{lemma:non-nullness-inequality} and assumption $(\mathbf{R2})$ to obtain the rough estimate 
\begin{align*}
    \abs{
        \hat{c}_{\Delta}(\eta^{z,i}, \xi_{\Delta})
        -
        \hat{c}_{\Delta}(\eta, \xi_{\Delta})
    }
    \leq 
    4 \frac{1}{\delta} e^R \sup_{\Delta \ni z, \xi_{\Delta}} \norm{c_{\Delta}(\cdot, \xi_\Delta)}_\infty. 
\end{align*}
In the case where $z$ is not contained in $\Delta$, we can be a bit more precise. Via the elementary algebraic rule 
\begin{align*}
    ac -bd = \frac{1}{2}[(a-b)(c+d) + (a+b)(c-d)],
\end{align*}
and Lemma~\ref{lemma:non-nullness-inequality} plus assumption $(\mathbf{R2})$ one obtains 
\begin{align*}
     \big\lvert c_{\Delta}(\xi_{\Delta}&\eta^{z,i}_{\Delta^c}, \eta_{\Delta})\big\rvert
        \abs{
            \frac{\gamma_{\Delta}(\xi_{\Delta}|\eta^{z,i}_{\Delta^c}) }{\gamma_{\Delta}(\eta_{\Delta}|\eta^{z,i}_{\Delta^c})}
            -
            \frac{\gamma_{\Delta}(\xi_{\Delta}|\eta_{\Delta^c}) }{\gamma_{\Delta}(\eta_{\Delta}|\eta_{\Delta^c})}
        }
        \\\
        +
        &\abs{\frac{\gamma_{\Delta}(\xi_{\Delta}|\eta_{\Delta^c}) }{\gamma_{\Delta}(\eta_{\Delta}|\eta_{\Delta^c})}}
        \abs{c_{\Delta}(\xi_{\Delta}\eta^{z,i}_{\Delta^c}, \eta_{\Delta})-c_{\Delta}(\xi_{\Delta}\eta_{\Delta^c}, \eta_{\Delta})}
    \\\
    =
    \frac{1}{2}
    \big\lvert c_{\Delta}&(\xi_{\Delta}\eta^{z,i}_{\Delta^c}, \eta_{\Delta})\big\rvert
    \abs{\frac{1}{\gamma_{\Delta}(\eta_{\Delta}|\eta^{z,i}_{\Delta^c})\gamma_{\Delta}(\eta_{\Delta}|\eta_{\Delta^c})}} 
     \abs{\gamma_{\Delta}(\xi_{\Delta}|\eta^{z,i}_{\Delta^c}) - \gamma_{\Delta}(\xi_{\Delta}|\eta_{\Delta^c})}
     \\\
     \times
    &\abs{\gamma_{\Delta}(\eta_{\Delta}|\eta^{z,i}_{\Delta^c})+\gamma_{\Delta}(\eta_{\Delta}|\eta_{\Delta^c})}
    +
    \abs{\frac{\gamma_{\Delta}(\xi_{\Delta}|\eta_{\Delta^c}) }{\gamma_{\Delta}(\eta_{\Delta}|\eta_{\Delta^c})}}
    \abs{c_{\Delta}(\xi_{\Delta}\eta^{z,i}_{\Delta^c}, \eta_{\Delta})-c_{\Delta}(\xi_{\Delta}\eta_{\Delta^c}, \eta_{\Delta})}
    \\\
    \leq 
    \frac{1}{\delta^2}&e^{2R}\norm{\nabla_z^i \gamma_\Delta (\cdot | \cdot)}_\infty 
    +
    \frac{1}{\delta}e^R\norm{\nabla_z^i c_\Delta(\cdot, \eta_\Delta)}_\infty. 
\end{align*}
Now assumptions $\mathbf{(R1)}-\mathbf{(R4)}$ and $\mathbf{(S1)}$ are sufficient to conclude that
\begin{align*}
        \sup_{y \in \Z^d}\sum_{\Delta \ni y}\sum_{z \neq y}\sum_{\xi_{\Delta}}\sum_{i=1}^q \norm{\nabla^i_z \hat{c}_\Delta(\cdot, \xi_\Delta)}_\infty < \infty. 
\end{align*}
\textit{Ad iii.: } 
This follows from the same estimates as in $ii.$ and  assumptions $\mathbf{(R3)}$ and $\mathbf{(S3)}$ by the equivalence that we stated before the proof. 
\newline 
\textit{Ad iv.: } 
This is clear by definition of the time-reversal rates $(\hat{c}_\Delta(\cdot, \xi_\Delta))_{\Delta \Subset \Z^d, \xi_\Delta \in \Omega_\Delta}$. 
\end{proof}

\begin{remark}
    The regularity statement in Lemma~\ref{lemma:regularity-time-reversal} in particular implies that the process with rates $\hat{c}_\Delta(\eta, \xi_\Delta)$ is well-defined. By using the switching Lemma~\ref{generalized-switching-lemma} one can then easily show that this process, with semigroup $(\hat{P}_t)_{t \geq 0}$, is dual to the original process in the sense that 
    for all $f,g \in C(\Omega)$ it holds that
        \begin{align*}
            \int_\Omega \left(P(t)f(\eta)\right) g(\eta)\mu(d\eta) = \int_\Omega f(\eta)(\hat{P}(t)g(\eta))\mu(d\eta), \quad t \geq 0. 
    \end{align*}
    Duality of Markov processes plays a big role in contemporary probability theory. In that context, the duality we have here is known as \textit{duality with respect to a measure}, see~\cite[Definition 1.3.]{jansen_notions_2014} and the remarks thereafter.  
    Studying the dual process can often yield useful information about the original process that is hard to obtain in other ways. However, at this point, we won't dive too deep into the analysis of the time-reversed interacting particle system and mainly use the time-reversal rates for notational simplicity.
    Investigating what can be done by following the duality approach further could be a path for future research. 
\end{remark}

Equipped with these estimates, we are now almost ready to show the main result of this section. We just need two more technical helpers to make our life a little easier. 
The first one is concerned with the density of Gibbs measures under certain transformations. 
\begin{lemma}\label{lemma:transformed-gibbs-measure}
    For $\Delta \Subset \Z^d$ and $\zeta_{\Delta}, \xi_{\Delta} \in \Omega_{\Delta}$ we define a map
    \begin{align*}
        G_{\xi_{\Delta}}: \Omega \to [\xi_{\Delta}]:=\{\omega\in\Omega\colon \omega_\Delta=\xi_\Delta\}, \quad \eta \mapsto \xi_{\Delta}\eta_{\Delta^c}. 
    \end{align*}
    Then, for $\mu \in \mathscr{G}(\gamma)$ we have that $\mu-a.s.$
    \begin{align}\label{density-of-transformed-measure}
        \frac{d(\mu \circ G_{\xi_{\Delta}}^{-1})}{d\mu}(\eta) = \mathbf{1}_{[\xi_{\Delta}]}(\eta)\sum_{\zeta_{\Delta}} \frac{\gamma_{\Delta}(\zeta_{\Delta} | \eta_{\Delta^c})}{\gamma_{\Delta}(\xi_{\Delta} | \eta_{\Delta^c})}, \quad \eta \in \Omega. 
    \end{align}
\end{lemma}

\begin{proof}
    Let $\Lambda \Subset \Z^d$ be such that $\Delta \subset \Lambda$ and let $\chi_{\Lambda} \in \Omega_{\Lambda}$. Then, we have 
\begin{align}\label{writing-out-conditional-expectation}
    (\mu \circ G_{\xi_{\Delta}}^{-1})([\xi_{\Delta}\chi_{\Lambda \setminus \Delta}])
    &=
    \sum_{\zeta_{\Delta}}\mu([\zeta_{\Delta}\chi_{\Lambda \setminus \Delta}])
    =
    \sum_{\zeta_{\Delta}}\int_{\Omega}\mathbf{1}_{\zeta_{\Delta}}(\omega)\mathbf{1}_{\chi_{\Lambda \setminus \Delta}}(\omega)\mu(d\omega) 
    \nonumber \\\
    &=
    \sum_{\zeta_{\Delta}}\int_{\Omega} \bbE[\mathbf{1}_{\zeta_{\Delta}} \mathbf{1}_{\chi_{\Lambda \setminus \Delta}} | \calF_{\Delta^c}](\omega) \mu(d\omega)
    \\\
    &= 
    \sum_{\zeta_{\Delta}}\int_{\Omega} \bbE[\mathbf{1}_{\zeta_{\Delta}} | \calF_{\Delta^c}](\omega) \mathbf{1}_{\chi_{\Lambda \setminus \Delta}}(\omega)\mu(d\omega).
    \nonumber 
\end{align} 
Now, since $\mu$ is a Gibbs measure with respect to the specification $\gamma$ we know that $\mu$-a.s 
\begin{align*}
    \bbE[\mathbf{1}_{\zeta_{\Delta}} | \calF_{\Delta^c}](\omega) = \gamma_{\Delta}(\zeta_{\Delta}|\omega_{\Delta^c}),
\end{align*}
and the right-hand-side is $\calF_{\Delta^c}$-measurable. 
Therefore, we can write each summand of~\eqref{writing-out-conditional-expectation} as 
\begin{align*}
\int_{\Omega} \bbE&[\mathbf{1}_{\zeta_{\Delta}} | \calF_{\Delta^c}](\omega) \mathbf{1}_{\chi_{\Lambda \setminus \Delta}}(\omega)\mu(d\omega)
    = \int_{\Omega} \gamma_{\Delta}(\zeta_{\Delta}|\omega_{\Delta^c}) \mathbf{1}_{\chi_{\Lambda \setminus \Delta}}(\omega)\mu(d\omega)
    \nonumber \\\
    &= 
    \int_{\Omega} \frac{\gamma_{\Delta}(\xi_{\Delta}| \omega_{\Delta^c})}{\gamma_{\Delta}(\xi_{\Delta}| \omega_{\Delta^c})}\gamma_{\Delta}(\zeta_{\Delta}|\omega_{\Delta^c}) \chi_{\Lambda \setminus \Delta}(\omega)\mu(d\omega)
    \\\
    &=
    \int_{\Omega} 
        \frac{\gamma_{\Delta}(\zeta_{\Delta}|\omega_{\Delta^c})}{\gamma_{\Delta}(\xi_{\Delta}| \omega_{\Delta^c})}
        \gamma_{\Delta}(\xi_{\Delta}| \omega_{\Delta^c}) 
        \chi_{\Lambda \setminus \Delta}(\omega)\mu(d\omega)
    \nonumber \\\
    &=\int_{\Omega}
        \bbE\left[ 
            \frac{\gamma_{\Delta}(\zeta_{\Delta}|\cdot_{\Delta^c})}{\gamma_{\Delta}(\xi_{\Delta}| \cdot_{\Delta^c})}
            \mathbf{1}_{\xi_{\Delta}}
            \mathbf{1}_{\chi_{\Lambda \setminus \Delta}}
        \Big| \calF_{\Delta^c}\right](\omega)
        \mu(d\omega)
    = \int_{\Omega}
            \frac{\gamma_{\Delta}(\zeta_{\Delta}|\omega_{\Delta^c})}{\gamma_{\Delta}(\xi_{\Delta}| \omega_{\Delta^c})}
            \mathbf{1}_{\xi_{\Delta}\chi_{\Lambda \setminus \Delta}}(\omega)
            \mu(d\omega). 
    \nonumber 
\end{align*}
Summing up over $\zeta_{\Delta}$ now gives us the claimed density. 
\end{proof}

The second technical result,
reminiscent of Lebesgue's differentiation theorem, will not only be used in this section but also in other parts of the manuscript. 
\begin{lemma}[Differentiation lemma]\label{differentiation-lemma}
    Let $\mu$ be a probability measure on $\Omega$ such that we have $\mu(\eta_{\Lambda})>0$ for all $\Lambda \Subset \Z^d$ and $\eta \in \Omega$. 
    Then, for any continuous functions $f: \Omega \to \R$ we have that for all $\eta \in \Omega$
    \begin{align*}
         \lim_{\Lambda \uparrow \Z^d} \frac{1}{\mu(\eta_{\Lambda})}\int\mathbf{1}_{\eta_{\Lambda}}(\xi)f(\xi)\mu(d\xi) = f(\eta). 
    \end{align*}
    Moreover, if $f$ is uniformly continuous, then the claimed convergence is also uniform in $\eta \in \Omega$. 
\end{lemma}

\begin{proof}
    First note that for fixed $\Lambda \Subset \Z^d$ we have the trivial inequalities 
    \begin{align}\label{sandwich}
        - \infty < \inf_{\xi: \xi_{\Lambda}=\eta_{\Lambda}}f(\xi) \leq f(\eta) \leq \sup_{\xi: \xi_{\Lambda}=\eta_{\Lambda}}f(\xi) < \infty. 
    \end{align}
    The continuity of $f$ implies that 
    \begin{align*}
        \lim_{\Lambda \uparrow \Z^d}\inf_{\xi: \xi_{\Lambda}=\eta_{\Lambda}}f(\xi) = f(\eta), \quad \lim_{\Lambda \uparrow \Z^d}\sup_{\xi: \xi_{\Lambda}=\eta_{\Lambda}}f(\xi) = f(\eta).
    \end{align*} 
    Combining this with~\eqref{sandwich} and the squeeze theorem (for nets) from real analysis yields
    \begin{align*}
        \lim_{\Lambda \uparrow \Z^d} \frac{1}{\mu(\eta_{\Lambda})}\int\mathbf{1}_{\eta_{\Lambda}}(\xi)f(\xi)\mu(d\xi) = f(\eta). 
    \end{align*}
    This concludes the proof. 
\end{proof}

\begin{proposition}\label{proposition:oscillation-equation}
    Assume that the rates $(c_{\Delta}(\cdot, \xi_{\Delta}))_{\Delta \Subset \Z^d, \xi_{\Delta}\in \Omega_{\Delta}}$ of an interacting particle system satisfy the conditions $\mathbf{(R1)}-\mathbf{(R4)}$ and that there exists a time-stationary measure $\mu$ such that $\mu \in \mathscr{G}(\gamma)$, where the specification $\gamma = (\gamma_{\Delta})_{\Delta \Subset \Z^d}$ satisfies the conditions $\mathbf{(S1)}-\mathbf{(S4)}$.
    Then, for all $\Lambda \Subset \Z^d$ and all $\eta \in \Omega$ it holds that
    \begin{align}\label{equilibrium-oscillations-1}
         \sum_{\Delta \Subset \Z^d}\sum_{\xi_{\Delta}}\nabla_{\Lambda}\left(c_{\Delta}(\cdot, \xi_{\Delta})- \hat{c}(\cdot, \xi_{\Delta})\right) (\eta) = 0, 
    \end{align}      
    where as before
    \begin{align*}
        \hat{c}_{\Delta}(\eta, \xi_{\Delta}) := c_{\Delta}(\xi_{\Delta}\eta_{\Delta^c}, \eta_{\Delta})\frac{\gamma_{\Delta}(\xi_{\Delta}|\eta_{\Delta^c})}{\gamma_{\Delta}(\eta_{\Delta}|\eta_{\Delta^c})}. 
    \end{align*}
    Moreover, we even have 
    \begin{align}\label{oscillation-sum-equation}
        \forall z \in \Z^d \ \forall i=1,\dots,q \  \forall \eta \in \Omega: \quad
        \sum_{\Delta \Subset \Z^d}\sum_{\xi_{\Delta}}\nabla_z^i\left(c_{\Delta}(\cdot, \xi_{\Delta})- \hat{c}(\cdot, \xi_{\Delta})\right) (\eta) = 0.
    \end{align}
\end{proposition}

\begin{proof}
    For fixed $\Lambda \Subset \Z^d$ define functions $h_{\Delta} : \Omega \to \R$ for each $\Delta \Subset \Z^d$ by 
    \begin{align*}
        h_{\Delta}(\eta) := \sum_{\zeta_{\Delta}}\nabla_{\Lambda} (c_{\Delta}(\cdot, \zeta_{\Delta})- \hat{c}_{\Delta}(\cdot, \zeta_{\Delta}))(\eta), \quad \eta \in \Omega.
    \end{align*}
    We have to show that $\sum_{\Delta \Subset \Z^d}h_\Delta = 0$. To do this, we will first show that 
    \begin{align*}
        \int_\Omega \left(\sum_{\Delta \Subset \Z^d}h_\Delta(\eta) \right)^2 \mu(d\eta) = 0
    \end{align*}
    and then use Lemma \ref{differentiation-lemma} to conclude that the integrand vanishes everywhere (and not just $\mu$-almost everywhere). 
    We start out by calculating $\int_{\Omega}h_\Delta h_\Theta \ d\mu$ for an arbitrary $\Theta \Subset \Z^d$. To do this, we will split the summation in two parts, namely 
    \begin{align*}
        \int_\Omega h_\Delta(\eta)h_\Theta(\eta) \ \mu(d\eta)
        =
        \sum_{\zeta_\Lambda}&\sum_{\xi_\Theta}
        \int_\Omega h_\Delta(\eta)
        \left(
        c_\Theta(\zeta_\Lambda \eta_{\Lambda^c}, \xi_\Theta) 
        -
        \hat{c}_\Theta(\zeta_\Lambda \eta_{\Lambda^c}, \xi_\Theta) 
        \right) 
        \mu(d\eta) 
        \\\
        -
        &\sum_{\zeta_\Lambda}\sum_{\xi_\Theta}
        \int_\Omega h_\Delta(\eta)
        \left(
            c_\Theta(\eta, \xi_\Theta)
            -
            \hat{c}_\Theta(\eta, \xi_\Theta)
        \right) 
        \mu(d\eta)
        =: \mathbf{I} + \mathbf{II}. 
    \end{align*}
   Using~\eqref{generalized-switching-lemma-identity} we can write the summands in $\mathbf{II}$ as 
    \begin{align*}
        \sum_{\xi_{\Theta}}
            \int_{\Omega} h_{\Delta}(\eta)(\hat{c}_{\Theta}(\eta, \xi_{\Theta})-c_{\Theta}(\eta, \xi_{\Theta}))\mu(d\eta)
         = 
         \sum_{\xi_{\Theta}}
            \int_{\Omega}c_{\Theta}(\eta, \xi_{\Theta})\left(h_{\Delta}(\xi_{\Theta}\eta_{\Theta^c})- h_{\Delta}(\eta)\right)\mu(d\eta). 
    \end{align*}
    For $\mathbf{I}$ we first apply Lemma \ref{lemma:transformed-gibbs-measure} to obtain 
    \begin{align*}
        \sum_{\xi_{\Theta}}\sum_{\zeta_{\Lambda}} &
            \int_{\Omega}
                h_{\Delta}(\omega)
                (
                    c_{\Theta}(\zeta_{\Lambda} \omega_{\Lambda^c},\xi_{\Theta})
                    -
                    \hat{c}_{\Theta}(\zeta_{\Lambda} \omega_{\Lambda^c},\xi_{\Theta})
                )
            \mu(d\omega)
        \nonumber \\\
        = 
        &\sum_{\xi_{\Theta}}\sum_{\zeta_{\Lambda}} \sum_{\eta_{\Lambda}}
            \int_{[\eta_{\Lambda}]}
                h_{\Delta}(\omega)
                (
                    c_{\Theta}(\zeta_{\Lambda} \omega_{\Lambda^c},\xi_{\Theta})
                    -
                    \hat{c}_{\Theta}(\zeta_{\Lambda} \omega_{\Lambda^c},\xi_{\Theta})
                )
            \mu(d\omega)
        \\\
        = 
        &\sum_{\xi_{\Theta}}\sum_{\zeta_{\Lambda}} \sum_{\eta_{\Lambda}}
            \int_{[\zeta_{\Lambda}]}
            h_{\Delta}(\eta_{\Lambda}\omega_{\Lambda^c})
            (
                c_{\Theta}(\omega, \xi_{\Theta})
                -
                \hat{c}_{\Theta}(\omega,\xi_{\Theta})
            )
            \left( \mu \circ G_{\zeta_{\Lambda}}^{-1}\right)(d\omega)
        \nonumber \\\
        = 
        &\sum_{\xi_{\Theta}}\sum_{\zeta_{\Lambda}} \sum_{\eta_{\Lambda}}
            \int_{[\zeta_{\Lambda}]}
                h_{\Delta}(\eta_{\Lambda}\omega_{\Lambda^c})
                \frac{\gamma_{\Lambda}(\eta_{\Lambda}| \omega_{\Lambda^c})}{\gamma(\zeta_{\Lambda}|\omega_{\Lambda^c})}
                (
                    c_{\Theta}(\omega, \xi_{\Theta})
                    -
                    \hat{c}_{\Theta}(\omega,\xi_{\Theta})
                )
            \mu(d\omega). 
            \nonumber
    \end{align*}
    To this we can now apply $(\ref{generalized-switching-lemma-identity})$ to write
    \begin{align*}
        \sum_{\xi_{\Theta}}\sum_{\zeta_{\Lambda}} \sum_{\eta_{\Lambda}}&
            \int_{[\zeta_{\Lambda}]}
                h_{\Delta}(\eta_{\Lambda}\omega_{\Lambda^c})
                \frac{\gamma_{\Lambda}(\eta_{\Lambda}| \omega_{\Lambda^c})}{\gamma(\zeta_{\Lambda}|\omega_{\Lambda^c})}
                (
                    c_{\Theta}(\omega, \xi_{\Theta})
                    -
                    \hat{c}_{\Theta}(\omega,\xi_{\Theta})
                )
            \mu(d\omega)
        \nonumber \\\
        = 
        &\sum_{\xi_{\Theta}}\sum_{\zeta_{\Lambda}} \sum_{\eta_{\Lambda}}
            \int_{\Omega}
                c_{\Theta}(\omega, \xi_{\Theta})
                \Big[
                    \mathbf{1}_{\zeta_{\Lambda}}(\omega)
                    h_{\Delta}(\eta_{\Lambda}\omega_{\Lambda^c})
                    \frac{\gamma_{\Lambda}(\eta_{\Lambda}| \omega_{\Lambda^c})}{\gamma(\zeta_{\Lambda}|\omega_{\Lambda^c})}
                    \\\
                    &\hspace{1cm}
                    -\mathbf{1}_{\zeta_{\Lambda}}(\xi_{\Theta}\omega_{\Theta^c})
                    h_{\Delta}(\xi_{\Theta \setminus \Lambda}\eta_{\Lambda}\omega_{\Lambda^c})
                    \frac{\gamma_{\Lambda}(\eta_{\Lambda}| \xi_{\Theta \setminus \Lambda}\eta_{\Lambda}\omega_{\Lambda^c})}{\gamma(\zeta_{\Lambda}|\xi_{\Theta \setminus \Lambda}\eta_{\Lambda}\omega_{\Lambda^c})}
                \Big]
            \mu(d\omega). 
    \end{align*}
    So, if we define a function $g_\Delta : \Omega \to \R$ by  
    \begin{align*}
        g_{\Delta}(\omega)
        :=
        \sum_{\zeta_{\Lambda}}\sum_{\eta_{\Lambda}}
            \mathbf{1}_{\zeta_{\Lambda}}(\omega)
            h_{\Delta}(\eta_{\Lambda}\omega_{\Lambda^c})
            \frac{\gamma_{\Lambda}(\eta_{\Lambda}| \omega_{\Lambda^c})}{\gamma_{\Lambda}(\zeta_{\Lambda}|\omega_{\Lambda^c})}, \quad \omega \in \Omega, 
    \end{align*}
    then combining the above calculations with the assumption that $\mu$ is time-stationary with respect to the Markovian dynamics generated by $\mathscr{L}$, implies that for all $\Delta \Subset \Z^d$ we have 
    \begin{align}\label{proof:oscillation-equation-1}
        \int_{\Omega}h_{\Delta}(\eta)\left(\sum_{\Theta \Subset \Z^d}h_{\Theta}(\eta)\right)\mu(d\omega)
        =
        \int_{\Omega}\mathscr{L}(h_{\Delta}-g_{\Delta})(\eta)\mu(d\eta)
        =
        0.
    \end{align}
    At this point, note that our assumptions on the rates and the specification ensure that 
    \begin{align*}
        h_{\Delta}, g_{\Delta} \in D(\Omega) \subset \text{dom}(\mathscr{L}),
    \end{align*}
    so applying the generator $\mathscr{L}$ to the function $(h_\Delta - g_\Delta)$ is a well-defined operation. 
    By summing~\eqref{proof:oscillation-equation-1} over all $\Delta \Subset \Z^d$ we obtain 
    \begin{align}\label{proof:oscillation-equation-2}
        \int_{\Omega}\left(\sum_{\Theta \Subset \Z^d}h_{\Theta}(\omega)\right)^2\mu(d\omega) = 0.
    \end{align}
    This only tells us that the desired equality~\eqref{equilibrium-oscillations-1} holds \textit{for $\mu$-a.e.~$\eta \in \Omega$}, which is not enough, as we will see later. However, the situation is not as dire as it may seem at first.
    Since we assumed that $\mu$ is non-null, we can use Lemma~\ref{differentiation-lemma} to conclude that~\eqref{equilibrium-oscillations-1} holds \textit{for every $\eta \in \Omega$}. Indeed, fix a sequence $(\Lambda_n)_{n \in \N}$ of finite subvolumes such that $\Lambda_n \uparrow \Z^d$. By non-nullness of $\mu$ and~\eqref{proof:oscillation-equation-2} we have for all $\eta \in \Omega$ and $n \in \N$
    \begin{align*}
        \frac{1}{\mu(\eta_{\Lambda_n})}\int \mathbf{1}_{\eta_{\Lambda_n}}(\omega)\left(\sum_{\Theta \Subset \Z^d}h_{\Theta}(\omega)\right)^2 \mu(d\omega) = 0. 
    \end{align*}
    Now note that our assumptions imply that the integrand in~\eqref{proof:oscillation-equation-2} is a continuous function. 
    By letting $n$ go to infinity and applying Lemma~\ref{differentiation-lemma} we see that for all $\omega \in \Omega$ we have 
    \begin{align*}
        \sum_{\Theta \Subset \Z^d}h_{\Theta}(\omega) = 0. 
    \end{align*}
    To see that we also have~\eqref{oscillation-sum-equation}, it suffices to note that
    for fixed $z \in \Z^d$ and $i\in \set{1,...,q}$ we can write
    \begin{align*}
        0 
        &= \sum_{\Delta \Subset \Z^d}\sum_{\xi_{\Delta}}
            \nabla_{\Lambda}\left(c_{\Delta}(\cdot, \xi_{\Delta})- \hat{c}(\cdot, \xi_{\Delta})\right) (\eta) 
        - \sum_{\Delta \Subset \Z^d}\sum_{\xi_{\Delta}}
            \nabla_{\Lambda}\left(c_{\Delta}(\cdot, \xi_{\Delta})- \hat{c}(\cdot, \xi_{\Delta})\right) (\eta^{z,i}) \nonumber \\\
        &= \sum_{\Delta \Subset \Z^d}\sum_{\xi_{\Delta}}\sum_{\zeta_{\Lambda}}\nabla_z^i\left(c_{\Delta}(\cdot, \xi_{\Delta})- \hat{c}(\cdot, \xi_{\Delta})\right) (\eta) \\\
        &= q^{\abs{\Lambda}}\sum_{\Delta \Subset \Z^d}\sum_{\xi_{\Delta}}\nabla_z^i\left(c_{\Delta}(\cdot, \xi_{\Delta})- \hat{c}(\cdot, \xi_{\Delta})\right) (\eta). \nonumber 
    \end{align*}
    This concludes the proof.
\end{proof}

\begin{proof}[Proof of Theorem \ref{theorem:all-gibbs-invariant}]
For $n \in \N$ we define the function 
    \begin{align*}
        \varphi_n(\eta) 
        := 
        \sum_{\Delta \cap \Delta_n \neq \emptyset}\sum_{\xi_{\Delta}}
            \left(
                c_\Delta(\eta, \xi_\Delta) - \hat{c}_\Delta(\eta, \xi_\Delta)
            \right),
    \end{align*}
    where $\Delta_n := [-n,n]^d$. Then, by Lemma~\ref{generalized-switching-lemma}, we have 
    \begin{align}\label{proof:all-gibbs-invariant-1}
        \int_\Omega \varphi_n(\eta)\mu(d\eta) 
        =
        0. 
    \end{align}
    For $\calF_\Lambda$-measurable $f:\Omega \to \R$ we can again use Lemma~\ref{generalized-switching-lemma} to get  
    \begin{align*}
        \int_\Omega \mathscr{L}f(\eta)\mu(d\eta)
        =
        \int_{\Omega}\sum_{\Delta \cap \Delta_n \neq \emptyset}\sum_{\xi_\Delta}c_\Delta(\eta, \xi_{\Delta})[f(\xi_\Delta \eta_{\Delta^c})-f(\eta)]
        =
        -\int_\Omega f(\eta)\varphi_n(\eta)\mu(d\eta),
    \end{align*}
    for all $n \in \N$ sufficiently large such that $\Lambda \subset \Delta_n$. Because of~\eqref{proof:all-gibbs-invariant-1} we can interpret the right-hand side as the (negative) covariance between $f$ and $\varphi_n$ under the probability measure $\mu$. 
    It remains to show that this covariance vanishes as $n$ tends to infinity. If $\varphi_n$ was $\calF_{\Delta_n^c}$-measurable, then we could directly use the standard covariance estimate \cite[Lemma 8.3.6] {durrett_probability_2019} to conclude this. 
    But in general this is not the case, and we need to proceed a bit more carefully. 
    First, observe that we can telescope $\varphi_n$ in the following elementary way
    \begin{align*}
        \varphi_n(\eta) 
        =
        \varphi_n(\mathbf{1}) + \sum_{m =0}^\infty ( \varphi_n(r_m \eta)-\varphi_n(r_{m+1}\eta)),
    \end{align*}
    where $\mathbf{1}\in \Omega$ is the configuration that is equal to $1$ at every site and the $\calF_{\Delta_{m-1}^c}$-measurable map $r_m :\Omega \to \Omega$ is defined by 
    \begin{align*}
        (r_m \eta)_x 
        = 
        \begin{cases}
            1, \quad &\text{if }\abs{x}_{\infty} < m,
            \\\
            \eta_x &\text{if } \abs{x}_{\infty} \geq m. 
        \end{cases} 
    \end{align*}
    Note that this construction also gives us 
    \begin{align*}
        0
        =
        \int_\Omega \varphi_n(\eta) \mu(d\eta)
        =  \sum_{m=0}^\infty \int_\Omega (\varphi_n(r_m \eta) - \varphi_n(r_{m+1}\eta))\mu(d\eta) + \varphi_n(\mathbf{1}),
    \end{align*}
    and therefore 
    \begin{align*}
        \int_\Omega f(\eta) \mu(d\eta) \cdot \left( \sum_{m=0}^\infty \int_\Omega (\varphi_n(r_m\eta) - \varphi_n(r_{m+1}\eta))\mu(d\eta) + \varphi_n(\mathbf{1}) \right)
        =
        0.
    \end{align*}
    This reduces our problem to estimating the covariance of $f$ and $(\varphi_n(r_m \cdot) - \varphi_n(r_{m+1}\cdot))$. 
    Observe that $(\varphi_n(r_m \cdot) - \varphi_n(r_{m+1}\cdot))$ is $\calF_{\Delta_{m-1}^c}$-measurable, so we can apply the classical estimate \cite[Lemma 8.3.6] {durrett_probability_2019} to see that
    \begin{align*}
        \abs{\int_{\Omega}f(\eta)\varphi_n(\eta)\mu(d\eta)}
        &=
        \abs{\int_{\Omega}f(\eta)\left(\varphi_n(\mathbf{1}) + \sum_{m =0}^\infty (\varphi_n(r_m \eta)-\varphi_n(r_{m+1}\eta))\right) \mu(d\eta)}
        \\\
        &\leq 
        4
        \sum_{m=0}^\infty 
            \norm{f}_\infty
            \norm{\varphi_n(r_m \cdot)-\varphi_n(r_{m+1}\cdot))}_\infty
            \alpha_{\mu}(\Lambda, m-1). 
    \end{align*}
    where we set $\alpha_\mu(\Lambda, -1) := 1$. 
    To estimate the summands of the series on the right-hand side, we first note that 
    \begin{align*}
        \norm{\varphi_n(r_m \cdot)-\varphi_n(r_{m+1}\cdot)}_\infty \leq \sum_{\abs{z} = m}\sum_{i=1}^q \norm{\nabla^i_z \varphi_n}_\infty. 
    \end{align*}
    This sum has $q \cdot O(m^{d-1})$ terms, because we only sum over the boundary sites of a $d$-dimensional hypercube. 
    For fixed $z\in \Z^d$ and $i \in \set{1, \dots, q}$ we can estimate the corresponding summand by using~\eqref{oscillation-sum-equation} to obtain
    \begin{align*}
        \norm{\nabla^i_z \varphi_n}_\infty
        &=
        \norm{\sum_{\Delta \cap \Delta_n \neq \emptyset}\sum_{\xi_\Delta}\nabla^i_z(c_\Delta(\cdot, \xi_\Delta)- \hat{c}_\Delta(\cdot, \xi_\Delta))}_\infty
        \\\
        &\leq
        \sum_{\Delta \cap \Delta_n = \emptyset}\sum_{\xi_\Delta}(\norm{\nabla^i_z c_\Delta(\cdot, \xi_\Delta)}_\infty + \norm{\nabla^i_z \hat{c}_\Delta(\cdot, \xi_\Delta)}_\infty)
        \leq 
        \beta(n). 
    \end{align*}
    In conjunction with the previous estimates, this gives us
    \begin{align*}
        \abs{\int_{\Omega}\mathscr{L}f(\eta)}
        &= 
        \abs{\int_\Omega f(\eta)\varphi_n(\eta)\mu(d\eta)}
        \\\
        &\leq 
        4\norm{f}_\infty \sum_{m=0}^{\infty}\alpha_\mu(\Lambda, m-1)\sum_{\abs{z}_\infty = m}\sum_{i=1}^q \norm{\nabla^i_z \varphi_n}_\infty
        \\\
        &\leq
        4 C \norm{f}_\infty \sum_{m=0}^{\infty}\alpha_\mu(\Lambda, m-1) m^{d-1} \beta(n). 
    \end{align*}
    By assumption, the term on the right-hand side vanishes as $n$ tends to infinity.
    Since this estimate holds for all local functions $f$, we can use dominated convergence to see that for all $g \in D(\Omega)$ it holds that
    \begin{align*}
        \int_\Omega \mathscr{L}g(\eta) \mu(d\eta) = 0. 
    \end{align*}
    But $D(\Omega)$ is a core for $\mathscr{L}$, so this implies that $\mu$ is time-stationary with respect to the Markovian dynamics generated by $\mathscr{L}$.
\end{proof}

Let us note that the condition on the mixing coefficients is hard to verify in practice and seems to be too strong, since we will see in the proof that we do not need to estimate the covariance of general pairs of functions with respect to $\mu$, but only for covariances of the form 
\begin{align*}
    \int_\Omega f(\eta) \varphi_n(\eta)\mu(d\eta).
\end{align*}

\subsection{Proof of Proposition~\ref{proposition:relative-entropy-loss-non-positive}}\label{section:proof-entropy-loss-non-positive}
We start out by deriving an explicit expression for the relative entropy loss in $\Lambda_n$ in terms of the rates and the measures $\mu, \nu$. 
\begin{lemma}\label{lemma:finite-volume-relative-entropy-loss}
    For $n\in \N$ and $\nu \in \calM_1(\Omega)$ we have 
    \begin{align*}
        g^n_{\mathscr{L}}(\nu|\mu) 
        =
        \sum_{\eta_{\Lambda_n}}\sum_{\Delta \cap \Lambda_n \neq \emptyset}\sum_{\xi_{\Delta}}
            \int \nu(d\omega)
                c_{\Delta}(\omega, \xi_{\Delta})
                \left[
                    \mathbf{1}_{\eta_{\Lambda_n}}(\xi_{\Delta}\omega_{\Delta^c}) 
                    -
                    \mathbf{1}_{\eta_{\Lambda_n}}(\omega)
                    \right]
            \log
            \left(
                \frac{\nu(\eta_{\Lambda_n})}{\mu(\eta_{\Lambda_n})}
            \right).
    \end{align*}
\end{lemma}
\begin{proof}
    This can be seen by a direct calculation using the definition of the generator. We have 
    \begin{align*}
        g^n_{\mathscr{L}}(\nu|\mu) 
        &=
        \sum_{\eta_{\Lambda_n}}
            \nu(\mathscr{L}\mathbf{1}_{\eta_{\Lambda_n}})
            \log
            \left(
                \frac{\nu(\eta_{\Lambda_n})}{\mu(\Lambda_n)}
            \right)
        \nonumber \\\
        &=
        \sum_{\eta_{\Lambda_n}}\sum_{\Delta \cap \Lambda_n \neq \emptyset}\sum_{\xi_{\Delta}}
        \int \nu(d\omega)
            c_{\Delta}(\omega, \xi_{\Delta})
            \left[
                \mathbf{1}_{\eta_{\Lambda_n}}(\xi_{\Delta}\omega_{\Delta^c}) 
                -
                \mathbf{1}_{\eta_{\Lambda_n}}(\omega)
                \right]
        \log
        \left(
            \frac{\nu(\eta_{\Lambda_n})}{\mu(\eta_{\Lambda_n})}
        \right).
    \end{align*}
    This completes the proof.
\end{proof}

To control the logarithmic terms in $g_\mathscr{L}^n(\nu \lvert \mu)$ 
we will again make use of Lemma \ref{lemma:non-nullness-inequality}. 
As already announced earlier, we will now take care of the terms corresponding to spin flips that happen outside $\Tilde{\Lambda}_n$. 
To do this properly, we will need the following simple combinatorial estimate. 

\begin{lemma}\label{combinatorial-lemma}
    Let $\Delta, \Lambda \Subset \Z^d$ be such that $\Delta \cap \Lambda \neq \emptyset$. Then, we have 
    \begin{align*}
        \abs{\set{x \in \Z^d: \ (\Delta + x)\cap \Lambda \neq \emptyset}} \leq \abs{\Delta}\abs{\Lambda}.
    \end{align*}
\end{lemma}

\begin{proof}
    Let $(\delta, \lambda) \in \Delta \times \Lambda$. Then, there exists a unique $x = x(\delta, \lambda) \in \Z^d$ such that
    \begin{align*}
        \delta + x = \lambda. 
    \end{align*}
    This clearly defines a surjective map 
    \begin{align*}
        \Delta \times \Lambda \ni (\delta, \lambda) \mapsto x(\delta, \lambda) \in \set{x \in \Z^d: \ (\Delta+x)\cap \Lambda \neq \emptyset}, 
    \end{align*}
    therefore we must have 
    \begin{align*}
        \abs{\set{x \in \Z^d: \ (\Delta + x)\cap \Lambda \neq \emptyset}} \leq \abs{\Delta}\abs{\Lambda},
    \end{align*}
as desired.
\end{proof}
\noindent 
This helps us in the following way. 
By assumption $\mathbf{(R2)}$, there are only finitely many different types of transitions, i.e., there are only finitely many distinct $\Delta \Subset \Z^d$ with $0 \in \Delta$ and $c_{\Delta}> 0$. Let $n$ be large enough such that all such basic shapes $\Delta$ are fully contained in $\Lambda_n$.
Then, the sum over all translations $(\Delta +x)$, $x \in \Z^d$,  of these basic shapes, such that $(\Delta +x) \cap \Lambda_n \neq \emptyset$ but $(\Delta + x) \nsubseteq \Tilde{\Lambda}_n$, has of the order $\abs{\Lambda_n \setminus \Tilde{\Lambda}_n}$ terms. 
In order to show that the boundary contributions are negligible in the density limit, it thus suffices to bound the terms uniformly. For this, we will again make use of Lemma \ref{lemma:non-nullness-inequality}. 

\begin{lemma}\label{lemma:boundary-contributions-boundary-order}
    Assume that the rates satisfy conditions $\mathbf{(R1)}$, $\mathbf{(R2)}$ and $\mathbf{(R4)}$. Moreover, assume that $\mu$ is time-stationary for the dynamics with $\mu \in \mathscr{G}(\gamma)$, such that the specification $\gamma$ satisfies $\mathbf{(S2)}$ and $\mathbf{(S4)}$. 
    For all $\nu \in \calM_1^{inv}(\Omega)$ and $n\in \N$ large enough, such that for all $\Delta \Subset \Z^d$ with $c_\Delta >0$ and $0 \in \Delta$ we have 
    $
        \Delta \Subset \Lambda_n,
    $
    it holds that 
    \begin{align*}
        \sum_{\eta_{\Lambda_n}}\sum_{\Delta \cap \Lambda_n \neq \emptyset: \Delta \nsubseteq \Tilde{\Lambda}_n} \sum_{\xi_{\Delta}}
            \int \nu(d\omega)
                c_{\Delta}(\omega, \xi_{\Delta})
                &\left[
                    \mathbf{1}_{\eta_{\Lambda_n}}(\xi_{\Delta}\omega_{\Delta^c}) 
                    -
                    \mathbf{1}_{\eta_{\Lambda_n}}(\omega)
                \right]
            \log
        \left(
            \frac{\nu(\eta_{\Lambda_n})}{\mu(\eta_{\Lambda_n})}
        \right)
        \\\
        &\leq 
        C \abs{\Lambda_n \setminus \Tilde{\Lambda}_n}
    \end{align*}
    for some constant $C > 0$ that does not depend on $n$ or $\nu$. 
    In particular, it holds that 
    \begin{align}\label{eqn:boundary-contributions}
        g_\mathscr{L}^n(\nu | \mu) \leq \Tilde{g}_\mathscr{L}^n(\nu | \mu) + o(\abs{\Lambda_n}). 
    \end{align}
    If $\nu$ is additionally non-null, then~\eqref{eqn:boundary-contributions} holds with equality.
\end{lemma}

\begin{proof}
    First assume that there is $\eta_{\Lambda_n} \in \Omega_{\Lambda_n}$ such that $\nu(\eta_{\Lambda_n})=0$ and $\nu(\xi_{\Delta\cap \Lambda_n}\eta_{\Lambda_n\setminus \Delta})>0$ for some $\xi_{\Delta}$.
    Then, the corresponding summand is equal to $-\infty$ and the upper bound is trivial. 
    If $\eta_{\Lambda_n}$ is such that $\nu(\eta_{\Lambda_n})=0$ and $\nu(\xi_{\Delta\cap \Lambda_n}\eta_{\Lambda_n\setminus \Delta})=0$ for all $\xi_{\Delta}$, then all the terms corresponding to $\eta_{\Lambda_n}$ are equal to $0$
    and we can therefore just omit them from the summation. 
    All in all, we can assume without loss of generality that $\nu(\eta_{\Lambda_n})>0$ for all $\eta_{\Lambda_n} \in \Omega_{\Lambda_n}$. 
    In this case, we can rearrange the sum we want to bound, without having to worry about adding and subtracting infinite terms, to obtain
    \begin{align*}
        \sum_{\eta_{\Lambda_n}}&\sum_{\Delta \cap \Lambda_n \neq \emptyset: \Delta \nsubseteq \Tilde{\Lambda}_n} \sum_{\xi_{\Delta}}
            \int \nu(d\omega)
                c_{\Delta}(\omega, \xi_{\Delta})
                \left[
                    \mathbf{1}_{\eta_{\Lambda_n}}(\xi_{\Delta}\omega_{\Delta^c}) 
                    -
                    \mathbf{1}_{\eta_{\Lambda_n}}(\omega)
                \right]
            \log
        \left(
            \frac{\nu(\eta_{\Lambda_n})}{\mu(\eta_{\Lambda_n})}
        \right)
        \nonumber \\\
        =
        &\sum_{\eta_{\Lambda_n}}\sum_{\Delta \cap \Lambda_n \neq \emptyset: \Delta \nsubseteq \Tilde{\Lambda}_n}\sum_{\zeta_{\Delta\setminus \Lambda_n}}\sum_{\psi_{\Delta \cap \Lambda_n}\neq \eta_{\Delta \cap \Lambda_n}}
            \int_{\Omega}\nu(d\omega)
                c_{\Delta}(\omega, \psi_{\Delta\cap \Lambda_n}\zeta_{\Delta \setminus \Lambda_n})
                \mathbf{1}_{\eta_{\Lambda_n}}(\omega) 
            \\\
            &\hspace{5cm}
            \times 
            \left[
                \log
                \left(
                    \frac{\nu(\psi_{\Delta \cap \Lambda_n}\eta_{\Lambda_n\setminus \Delta})}{\mu(\psi_{\Delta \cap \Lambda_n}\eta_{\Lambda_n \setminus\Delta})}
                \right)
                -
                \log
                \left(
                    \frac{\nu(\eta_{\Lambda_n})}{\mu(\eta_{\Lambda_n})}
                \right)
            \right]
        \nonumber \\\
        =
        &\sum_{\eta_{\Lambda_n}}\sum_{\Delta \cap \Lambda_n \neq \emptyset: \Delta \nsubseteq \Tilde{\Lambda}_n}\sum_{\zeta_{\Delta\setminus \Lambda_n}}\sum_{\psi_{\Delta \cap \Lambda_n}\neq \eta_{\Delta \cap \Lambda_n}}
            \int_{\Omega}\nu(d\omega)
                c_{\Delta}(\omega, \psi_{\Delta\cap \Lambda_n}\zeta_{\Delta \setminus \Lambda_n})
                \mathbf{1}_{\eta_{\Lambda_n}}(\omega) 
            \\\
            &\hspace{5cm}
            \times \log\left(
                \frac{\nu(\psi_{\Delta\cap\Lambda_n}\eta_{\Lambda_n \setminus \Delta})}{\nu(\eta_{\Lambda_n})}
            \right)
            \nonumber 
        \\\
        +
        &\sum_{\eta_{\Lambda_n}}\sum_{\Delta \cap \Lambda_n \neq \emptyset: \Delta \nsubseteq \Tilde{\Lambda}_n}\sum_{\zeta_{\Delta\setminus \Lambda_n}}\sum_{\psi_{\Delta \cap \Lambda_n}\neq \eta_{\Delta \cap \Lambda_n}}
            \int_{\Omega}\nu(d\omega)
                c_{\Delta}(\omega, \psi_{\Delta\cap \Lambda_n}\zeta_{\Delta \setminus \Lambda_n})
                \mathbf{1}_{\eta_{\Lambda_n}}(\omega) 
            \\\
            &\hspace{5cm}
            \times 
            \log\left(
                \frac{\mu(\eta_{\Lambda_n})}{\mu(\psi_{\Delta\cap\Lambda_n}\eta_{\Lambda_n \setminus \Delta})}
            \right)
        \nonumber 
        \\\
        =: &\, \mathbf{I} + \mathbf{II}. 
        \nonumber 
    \end{align*}
    We can now bound these two terms separately, starting with $\mathbf{II}$. Here we can apply Lemma~\ref{lemma:non-nullness-inequality} and use the translation invariance of the rates to get
    \begin{align*}
        \abs{\mathbf{II}} \leq \sup_{\Delta \ni 0: \xi_{\Delta}}\norm{c_{\Delta}(\cdot, \xi_{\Delta})}_{\infty}q^{\Delta} \log\left(\frac{1}{\delta(\mu)}\right) \abs{\Lambda_n \setminus \hat{\Lambda}_n} = o(\abs{\Lambda_n}),
    \end{align*}
    where $\delta(\mu)$ is the constant in the non-nullness estimate for $\mu$ and 
    \begin{align*}
        \hat{\Lambda}_n = [-2^n + n + L + 1, 2^n - n - L -1]^d, \quad  \text{with } L:=\max_{\Delta \ni 0: c_{\Delta}>0}\text{diam}(\Delta)+1.
    \end{align*}
    If $\nu$ is also non-null, then we can estimate $\mathbf{I}$ in exactly the same way. For general $\nu$, we first use the trivial estimate 
    \begin{align*}
        \log \leq \log^+ := \max\left\{0, \log(\cdot)\right\}
    \end{align*}
    and then 
    \begin{align*}
        x \log^+\left(\frac{1}{x}\right) \leq e^{-1}, \quad \forall x > 0,
    \end{align*}
    to get the upper bound 
    \begin{align*}
        \mathbf{I} 
        \leq 
        \sup_{\Delta \ni 0: \xi_{\Delta}}\norm{c_{\Delta}(\cdot, \xi_{\Delta})}_{\infty}
        &\sum_{\eta_{\Lambda_n}}\sum_{\Delta \cap \Lambda_n \neq \emptyset: \Delta \nsubseteq \Tilde{\Lambda}_n}\sum_{\zeta_{\Delta\setminus \Lambda_n}}\sum_{\psi_{\Delta \cap \Lambda_n}\neq \eta_{\Delta \cap \Lambda_n}}
            \nu({\eta_{\Lambda_n}})
            \log^+\left(
                \frac{\nu(\psi_{\Delta\cap\Lambda_n}\eta_{\Lambda_n \setminus \Delta})}{\nu(\eta_{\Lambda_n})}
            \right)
        \nonumber \\\
        =
        \sup_{\Delta \ni 0: \xi_{\Delta}}\norm{c_{\Delta}(\cdot, \xi_{\Delta})}_{\infty}
        &\sum_{\eta_{\Lambda_n}}\sum_{\Delta \cap \Lambda_n \neq \emptyset: \Delta \nsubseteq \Tilde{\Lambda}_n}\sum_{\zeta_{\Delta\setminus \Lambda_n}}\sum_{\psi_{\Delta \cap \Lambda_n}\neq \eta_{\Delta \cap \Lambda_n}}
            \nu(\psi_{\Delta\cap\Lambda_n}\eta_{\Lambda_n \setminus \Delta})
            \frac{\nu(\eta_{\Lambda_n})}{\nu(\psi_{\Delta\cap\Lambda_n}\eta_{\Lambda_n \setminus \Delta})}
            \nonumber 
            \\\
            &\hspace{6cm}
            \times 
            \log^+\left(
                \frac{\nu(\psi_{\Delta\cap\Lambda_n}\eta_{\Lambda_n \setminus \Delta})}{\nu(\eta_{\Lambda_n})}
            \right)
        \\\
        \leq 
        \sup_{\Delta \ni 0: \xi_{\Delta}}\norm{c_{\Delta}(\cdot, \xi_{\Delta})}_{\infty}
        & e^{-1}
        \abs{\Lambda_n \setminus \hat{\Lambda}_n} = o(\abs{\Lambda_n}). 
        \nonumber
    \end{align*}
    This completes the proof.
\end{proof}

Motivated by Lemma \ref{lemma:boundary-contributions-boundary-order}, we now define the \textit{approximating relative entropy loss in $\Lambda_n$} by 
\begin{align*}
    \Tilde{g}^n_{\mathscr{L}}(\nu|\mu)
    := 
    \sum_{\eta_{\Lambda_n}}\sum_{\Delta \subset \Tilde{\Lambda}_n} \sum_{\xi_{\Delta}}
            \int \nu(d\omega)
                c_{\Delta}(\omega, \xi_{\Delta})
                \left[
                    \mathbf{1}_{\eta_{\Lambda_n}}(\xi_{\Delta}\omega_{\Delta^c}) 
                    -
                    \mathbf{1}_{\eta_{\Lambda_n}}(\omega)
                \right]
            \log
        \left(
            \frac{\nu(\eta_{\Lambda_n})}{\mu(\eta_{\Lambda_n})}
        \right).
\end{align*}
\noindent 
This is the rate at which the relative entropy in $\Lambda_n$ changes due to spin flips inside the smaller region $\Tilde{\Lambda}_n \subset \Lambda_n$. 

\medskip
We now proceed by rewriting the approximating entropy loss in a way that makes it obvious that it is non-negative – up to some negligible terms that vanish when taking the density limit. Conceptually this representation is analogous to~\eqref{eqn:relative-entropy-loss-finite-state-spaces}  for finite state spaces, but for infinite systems we have to overcome some additional technical difficulties.  Some of these are already present in Holley's seminal work \cite{holley_free_1971}, but the non-reversibility forces us to work substantially harder.
Since we cannot apply the detailed-balance equations, we will have to rely on~\eqref{oscillation-sum-equation} to show that certain error terms are of boundary order. 

\medskip 
Before we start the proof, we first define  
\begin{align*}
    F_0(u) := 
    \begin{cases}
        u - u \log(u)-1, &\quad \text{if $u>0$}, \\\
        -1, &\quad  \text{otherwise},
    \end{cases}
\end{align*}
and
\begin{align*}
    F(\nu, n,\eta, \xi_{\Delta})
    := 
    \begin{cases}
        F_0\left(\frac{\nu(\eta_{\Lambda_n})}{\nu(\xi_{\Delta}\eta_{\Lambda_n \setminus \Delta})}\frac{\mu(\xi_{\Delta}\eta_{\Lambda_n \setminus \Delta})}{\mu(\eta_{\Lambda_n})}\right)\nu(\xi_{\Delta}\eta_{\Lambda_n \setminus \Delta}), \ \text{if $\nu(\xi_{\Delta}\eta_{\Lambda_n\setminus \Delta})>0$},\\\
        - \infty,  \hspace{4.5cm}\text{if $\nu(\xi_{\Delta}\eta_{\Lambda_n \setminus \Delta})=0$ and $\nu(\eta_{\Lambda_n})>0$},\\\
        0, \hspace{5cm}\text{if $\nu(\xi_{\Delta}\eta_{\Lambda_n\setminus \Delta})=\nu(\eta_{\Lambda_n})=0$.}
    \end{cases}
\end{align*}
Note that $F_0(\cdot)$ is non-positive, concave, and only vanishes at $u=1$. 
For a configuration $\omega \in \Omega_{\Lambda_n}$ (or $\omega \in \Omega$) let $r_n\omega$ denote the configuration defined by 
\begin{align}
    (r_n\omega)_y 
    =
    \begin{cases}
        \omega_y, \quad &\text{if } y \in \Lambda_n, \\\
        1, &\text{otherwise}.
    \end{cases}
\end{align}
This will serve as an infinite-volume extension of the finite-volume configuration $\omega$ (or a finite-volume approximation of the infinite-volume configuration $\omega$, depending on the point of view). To see that certain error terms are of boundary order, we will make use of the following lemma multiple times. 

\begin{lemma}\label{lemma:uniform-convergence-conditional-probabilities}
    Assume that the specification $\gamma$ satisfies conditions $\mathbf{(S1)}-\mathbf{(S2)}$.
    Let $\Delta \Subset \Z^d$ and fix $\xi_\Delta \in \Omega_\Delta$. Then, the following convergence holds uniform in $\eta \in \Omega$
    \begin{align*}
        \frac{\mu(\eta_{\Lambda_n})}{\mu(\xi_\Delta \eta_{\Lambda_n \setminus \Delta})} 
        \to 
        \frac{\gamma_\Delta(\eta_\Delta| \eta_{\Delta^c})}{\gamma_\Delta(\xi_\Delta| \eta_{\Delta^c})} 
        \quad 
        \text{as } n \to \infty.  
    \end{align*}
\end{lemma}

\begin{proof}
    As a first step, note that we can write
    \begin{align*}
        \frac{\mu(\eta_{\Lambda_n})}{\mu(\xi_\Delta \eta_{\Lambda_n \setminus \Delta})} 
        =
        \frac{\mu(\eta_\Delta| \eta_{\Lambda_n \setminus \Delta})}{\mu(\xi_\Delta | \eta_{\Lambda_n \setminus \Delta})}. 
    \end{align*}
    We first show that both the denominator and the numerator converge uniformly in $\eta$. For this, observe that the DLR equations imply 
    \begin{align*}
        \mu(\eta_\Delta| \eta_{\Lambda_n \setminus \Delta})
        &=
        \frac{1}{\mu(\eta_{\Lambda_n \setminus \Delta})}
        \int_{\Omega}\gamma_{\Delta}(\eta_{\Lambda_n}| \omega_{\Delta^c})\mu(d\omega)
        \\\
        &=
        \frac{1}{\mu(\eta_{\Lambda_n \setminus \Delta})}
        \int_{\Omega}\mathbf{1}_{\eta_{\Lambda_n \setminus \Delta}}(\omega)\gamma_\Delta(\eta_{\Delta}| \omega_{\Delta^c})\mu(d\omega).
    \end{align*}
    By Lemma~\ref{differentiation-lemma} and the uniform continuity of $\gamma_{\Delta}$ with respect to the boundary condition this implies the uniform convergence  
    \begin{align*}
        \mu(\eta_\Delta| \eta_{\Lambda_n \setminus \Delta}) 
        \to 
        \gamma_{\Delta}(\eta_{\Delta}| \eta_{\Delta^c})
        \quad
        \text{as } n \to \infty. 
    \end{align*}
    Exactly the same argument implies the uniform convergence of the denominator
    \begin{align*}
        \mu(\xi_\Delta| \eta_{\Lambda_n \setminus \Delta}) 
        \to 
        \gamma_{\Delta}(\xi_{\Delta}| \eta_{\Delta^c})
        \quad
        \text{as } n \to \infty. 
    \end{align*}
    Now we can again use the simple algebraic rule 
    \begin{align*}
        ad - bc = \frac{1}{2}[(a-b)(c+d)-(a+b)(c-d)]
    \end{align*}
    in conjunction with the non-nullness of $\gamma$, and hence $\mu$, to obtain the inequality
    \begin{align}\label{proof:uniform-convergence-conditional-distributions}
        &\abs{
            \frac{\mu(\eta_{\Lambda_n})}{\mu(\xi_\Delta \eta_{\Lambda_n \setminus \Delta})} 
            -
            \frac{\gamma_\Delta(\eta_\Delta| \eta_{\Delta^c})}{\gamma_\Delta(\xi_\Delta| \eta_{\Delta^c})} 
        }
        \\\
        &\leq \ 
        \frac{1}{\delta^2}e^{2 \abs{\Delta}}
        \left(
        \abs{\mu(\eta_\Delta | \eta_{\Lambda_n \setminus \Delta})-\gamma_\Delta(\eta_\Delta| \eta_\Delta^c)}
        + 
        \abs{\mu(\xi_\Delta | \eta_{\Lambda_n \setminus \Delta})-\gamma_\Delta(\xi_\Delta| \eta_\Delta^c)}
        \right). 
        \nonumber 
    \end{align}
    By our previous considerations, the right-hand side of~\eqref{proof:uniform-convergence-conditional-distributions} converges to zero uniformly in $\eta$ as $n$ tends to infinity. 
    The uniformity in $\Delta \Subset \Z^d$ and $\xi_\Delta$ such that $c_\Delta > 0$ is now a consequence of the assumption that there are only finitely many types of transitions. 
\end{proof}
\begin{lemma}\label{lemma:bulk-contribution}
    Assume that the rates satisfy conditions $\mathbf{(R1)}-\mathbf{(R5)}$ and that $\mu$ is time-stationary for the dynamics with $\mu \in \mathscr{G}(\gamma)$, and the specification $\gamma$ satisfies $\mathbf{(S1)}-\mathbf{(S3)}$. 
    Then for all $n \in \N$ and $\nu \in \calM_1(\Omega)$ we have 
    \begin{align}
        \Tilde{g}^n_{\mathscr{L}}(\nu|\mu)
        = 
        \sum_{\eta_{\Lambda_n}}\sum_{\Delta \subset \Tilde{\Lambda}_n}\sum_{\xi_{\Delta} \neq \eta_{\Delta}} 
        F(\nu, n, \eta_{\Lambda_n}, \xi_{\Delta})c^{(n)}_{\Delta}(\eta_{\Lambda_n},\xi_{\Delta})
        \frac{\mu(\eta_{\Lambda_n})}{\mu(\xi_{\Delta}\eta_{\Lambda_n \setminus \Delta})} + o(\abs{\Lambda_n}),
    \end{align}
    where we use the truncated rates 
    \begin{align}
        c^{(n)}_{\Delta}(\eta_{\Lambda_n}, \xi_{\Delta}) := 
        \begin{cases}
            \frac{1}{\nu(\eta_{\Lambda_n})}\int \mathbf{1}_{\eta_{\Lambda_n}}(\omega) c_{\Delta}(\omega, \xi_{\Delta})\nu(d\omega), \quad &\text{\normalfont if $\nu(\eta_{\Lambda_n})>0$,} \\\
            c_{\Delta}(r_n \eta_{\Lambda_n}, \xi_{\Delta}), &\text{\normalfont otherwise}.
        \end{cases}
    \end{align}
\end{lemma}

\begin{proof}
    If there is $\eta_{\Lambda_n} \in \Omega_{\Lambda_n}$ such that $\nu(\eta_{\Lambda_n}) = 0$ and $\nu(\xi_{\Delta}\eta_{\Lambda_n \setminus \Delta})>0$
     for some $\Delta \subset \Tilde{\Lambda}_n$ and $\xi_{\Delta}\neq \eta_{\Delta}$, 
     then equality holds in the sense that $-\infty = - \infty$. 
    If $\eta$ is such that $\nu(\eta_{\Lambda_n}) = 0$ and $\nu(\xi_{\Delta}\eta_{\Lambda_n \setminus \Delta})=0$ for all $\Delta, \xi_{\Delta}$, 
    then all of the corresponding terms are equal to $0$ on both sides with the convention $0 \log 0 = 0$.  
    Therefore, we can assume without loss of generality that $\nu(\eta_{\Lambda_n})>0$ for all $\eta \in \Omega_{\Lambda_n}$.
    This allows us to express $\Tilde{g}_\mathscr{L}^n(\nu| \mu)$ as  
    \begin{align}\label{proof:rewriting-bulk-contribution-1}
        &\sum_{\eta_{\Lambda_n}} \sum_{\Delta \subset \Tilde{\Lambda}_n} \sum_{\xi_{\Delta} \neq \eta_{\Delta}}
        \left[
            \int\mathbf{1}_{\xi_{\Delta}\eta_{\Lambda_n \setminus \Delta}}(\omega)c_{\Delta}(\omega, \eta_{\Delta})\nu(d\omega) 
            - \int\mathbf{1}_{\eta_{\Lambda_n}}(\omega)c_{\Delta}(\omega, \xi_{\Delta})\nu(d\omega) 
        \right]
        \nonumber
        \\\
        &\hspace{7cm}
        \times 
        \log\left(
            \frac{\nu(\eta_{\Lambda_n})}{\mu(\eta_{\Lambda_n})}
        \right) \nonumber 
        \\\
        &= \sum_{\eta_{\Lambda_n}} \sum_{\Delta \subset \Tilde{\Lambda}_n} \sum_{\xi_{\Delta} \neq \eta_{\Delta}}
            \int \mathbf{1}_{\eta_{\Lambda_n}}(\omega)c_{\Delta}(\omega, \xi_{\Delta})\nu(d\omega)
            \left[
                \log\left(
                    \frac{\nu(\xi_{\Delta}\eta_{\Lambda_n \setminus \Delta})}{\mu(\xi_{\Delta}\eta_{\Lambda_n \setminus \Delta})}
                \right)
                -
                \log\left(
                    \frac{\nu(\eta_{\Lambda_n})}{\mu(\eta_{\Lambda_n})}
                \right)
            \right]
        \nonumber \\\
        &= \sum_{\eta_{\Lambda_n}} \sum_{\Delta \subset \Tilde{\Lambda}_n} \sum_{\xi_{\Delta} \neq \eta_{\Delta}}
            c^{(n)}_{\Delta}(\eta_{\Lambda_n}, \xi_{\Delta}) \nu(\eta_{\Lambda_n})
            \log\left(
                \frac{\nu(\xi_{\Delta}\eta_{\Lambda_n \setminus \Delta})}{\mu(\xi_{\Delta}\eta_{\Lambda_n \setminus \Delta})}
                \frac{\mu(\eta_{\Lambda_n})}{\nu(\eta_{\Lambda_n})}
            \right)
        \\\
        &= \sum_{\eta_{\Lambda_n}} \sum_{\Delta \subset \Tilde{\Lambda}_n} \sum_{\xi_{\Delta} \neq \eta_{\Delta}}
            F(\nu, n, \eta_{\Lambda_n}, \xi_{\Delta}) c^{(n)}_{\Delta}(\eta_{\Lambda_n}, \xi_{\Delta}) 
            \frac{\mu(\eta_{\Lambda_n})}{\mu(\xi_{\Delta}\eta_{\Lambda_n \setminus \Delta})} 
        \nonumber \\\
        &\quad+ 
        \sum_{\eta_{\Lambda_n}} \sum_{\Delta \subset \Tilde{\Lambda}_n} \sum_{\xi_{\Delta} \neq \eta_{\Delta}}
        \Bigg[
            \int\mathbf{1}_{\eta_{\Lambda_n}}(\omega)c_{\Delta}(\omega, \xi_{\Delta})\nu(d\omega)
            \nonumber \\\
            &\hspace{4cm}
            -
            \frac{\mu(\eta_{\Lambda_n})}{\mu(\xi_{\Delta}\eta_{\Lambda_n \setminus \Delta})}
            \frac{\nu(\xi_{\Delta}\eta_{\Lambda_n \setminus \Delta})}{\nu(\eta_{\Lambda_n})}
            \int\mathbf{1}_{\eta_{\Lambda_n}}(\omega)c_{\Delta}(\omega, \xi_{\Delta})\nu(d\xi)
        \Bigg]. 
        \nonumber 
    \end{align}
    It remains to show that the second term on the right-hand side of~\eqref{proof:rewriting-bulk-contribution-1} is negligible when taking the density limit. 
    We do this by first decomposing the error term into three parts
    and estimating them all separately. More precisely, we write
    \begin{align*}
        \sum_{\eta_{\Lambda_n}} \sum_{\Delta \subset \Tilde{\Lambda}_n} \sum_{\xi_{\Delta} \neq \eta_{\Delta}}
        &\left[
            \int\mathbf{1}_{\eta_{\Lambda_n}}(\omega)c_{\Delta}(\omega, \xi_{\Delta})\nu(d\omega)
            - 
            \tfrac{\mu(\eta_{\Lambda_n})\nu(\xi_{\Delta}\eta_{\Lambda_n \setminus \Delta})
            }{\mu(\xi_{\Delta}\eta_{\Lambda_n \setminus \Delta})\nu(\eta_{\Lambda_n})}
            \int\mathbf{1}_{\eta_{\Lambda_n}}(\omega)c_{\Delta}(\omega, \xi_{\Delta})\nu(d\xi)
        \right]
        \\\
        &= \mathbf{I} + \mathbf{II} + \mathbf{III},
    \end{align*}
    where 
    \begin{align*}
        \mathbf{I}:&= 
        \sum_{\eta_{\Lambda_n}} \sum_{\Delta \subset \Tilde{\Lambda}_n} \sum_{\xi_{\Delta} \neq \eta_{\Delta}}
        \int\mathbf{1}_{\eta_{\Lambda_n}}(\omega)
        \left(
            c_{\Delta}(\omega, \xi_{\Delta}) - c_{\Delta}(r_n \omega, \xi_{\Delta})
        \right)\nu(d\omega),
        \\\
        \mathbf{II}:&= 
        \sum_{\eta_{\Lambda_n}} \sum_{\Delta \subset \Tilde{\Lambda}_n} \sum_{\xi_{\Delta} \neq \eta_{\Delta}}
        \tfrac{\nu(\xi_{\Delta}\eta_{\Lambda_n \setminus \Delta})}{\nu(\eta_{\Lambda_n})}
        \int\mathbf{1}_{\eta_{\Lambda_n}}(\omega)
        \left(
            \hat{c}_{\Delta}(r_n \omega, \xi_{\Delta})-c_{\Delta}(\omega, \xi_{\Delta})\tfrac{\mu(\eta_{\Lambda_n})}{\mu(\xi_{\Delta}\eta_{\Lambda_n \setminus \Delta})}
        \right)\nu(d\omega), 
        \\\
        \mathbf{III}:&=
        \sum_{\eta_{\Lambda_n}} \sum_{\Delta \subset \Tilde{\Lambda}_n} \sum_{\xi_{\Delta} \neq \eta_{\Delta}}
        \nu(\eta_{\Lambda_n})c_{\Delta}(r_n \eta, \xi_{\Delta})
        - 
        \sum_{\eta \in \Omega_{\Lambda_n}} \sum_{\Delta \subset \Tilde{\Lambda}_n} \sum_{\xi_{\Delta} \neq \eta_{\Delta}}
        \nu(\xi_{\Delta}\eta_{\Lambda_n \setminus \Delta})\hat{c}_{\Delta}(r_n \eta_{\Lambda_n}, \xi_{\Delta}).
    \end{align*}
    Here, we used that the function $r_n: \Omega \to \Omega$ is constant on the cylinder sets $\set{\omega: \ \omega_{\Lambda_n} = \eta_{\Lambda_n}}$ for $\eta_{\Lambda_n} \in \Omega_{\Lambda_n}$. We now estimate these three terms separately. 
    The term $\mathbf{I}$ is of the order $o(\abs{\Lambda_n})$, because the integrands are of the order $o(1)$ by uniform continuity of the rates. 
    Similarly, the term $\mathbf{II}$ is of the order $o(\abs{\Lambda_n})$ because of Lemma~\ref{lemma:uniform-convergence-conditional-probabilities} and the definition of the rates of the time-reversal. Note that the convergence in Lemma~\ref{lemma:uniform-convergence-conditional-probabilities} is uniform over all summands by assumptions $\mathbf{(R2)}$ and $\mathbf{(R4)}$. 
    The term $\mathbf{III}$ would be zero in the reversible case and needs some extra attention. First note that we can rewrite it as 
    \begin{align*}
        \mathbf{III}
        &=
        \sum_{\eta_{\Lambda_n}} 
        \nu(\eta_{\Lambda_n})
        \sum_{\Delta \subset \Tilde{\Lambda}_n} \sum_{\xi_{\Delta} \neq \eta_{\Delta}}
        \left(
            c_{\Delta}(r_n \eta_{\Lambda_n}, \xi_{\Delta})-\hat{c}_{\Delta}(r_n \eta_{\Lambda_n}, \xi_{\Delta})
        \right)
        \\\
        &=:
        \sum_{\eta_{\Lambda_n} \in \Omega_{\Lambda_n}}
        \nu(\eta_{\Lambda_n}) \varphi_n(\eta_{\Lambda_n}). 
    \end{align*}
    Now for fixed $\Delta \subset \Tilde{\Lambda}_n$ we have by definiton of the time-reversal rates $\hat{c}$
    \begin{align}
        \sum_{\eta_{\Lambda_n}} &\sum_{\xi_{\Delta} \neq \eta_{\Delta}}
        \gamma_{\Lambda_n}(\eta_{\Lambda_n}\lvert r_n \eta_{\Lambda_n^c})
        \left(
            c_{\Delta}(r_n \eta_{\Lambda_n}, \xi_{\Delta})
            -
            \hat{c}_{\Delta}(r_n \eta_{\Lambda_n}, \xi_{\Delta})
        \right)
        \nonumber \\\
        =
        &\sum_{\eta_{\Lambda_n}} \sum_{\xi_{\Delta} \neq \eta_{\Delta}}
        \Big(
            \gamma_{\Lambda_n}(\eta_{\Lambda_n}\lvert r_n \eta_{\Lambda_n^c})
            c_{\Delta}(r_n \eta, \xi_{\Delta})
            \\\
            &\hspace{2cm} -\gamma_{\Lambda_n}(\xi_{\Delta}\eta_{\Lambda_n\setminus \Delta}\lvert r_n \xi_{\Delta}\eta_{\Lambda_n\setminus\Delta^c})
            c_{\Delta}(\xi_{\Delta}r_n\eta_{\Delta^c}, \eta_{\Delta})
        \Big)
        = 
        0. \nonumber 
    \end{align}
    So by summing over $\Delta \subset \Tilde{\Lambda}_n$ we see that 
    \begin{align*}
        \sum_{\eta_{\Lambda_n}}\gamma_{\Lambda_n}(\eta_{\Lambda_n}|r_n \eta_{\Lambda_n^c})\varphi_n(\eta_{\Lambda_n}) = 0. 
    \end{align*}
    Since $\gamma_{\Lambda_n}(\eta_{\Lambda_n}| r_n\eta_{\Lambda_n^c})>0$ for all $\eta$, by non-nullness of $\gamma$, we can conclude that 
    \begin{align}
        \inf_{\omega}\varphi_n(\omega) \leq 0 \leq \sup_{\omega}\varphi_n(\omega). 
    \end{align}
    This allows us to estimate the supremum norm of $\varphi_n$ by its oscillations. This yields
    \begin{align}\label{proof:rewriting-bulk-contribution-2}
        \norm{\varphi_n}_{\infty} 
        &\leq
        \sum_{z \in \Lambda_n}\sum_{i=1}^q\norm{\nabla^i_z \varphi_n}_{\infty}
        \nonumber 
        \\\
        &=
        \sum_{z \in \Lambda_n} \sum_{i=1}^q
        \norm{
            \sum_{\Delta \subset \Tilde{\Lambda}_n} \sum_{\xi_{\Delta}} 
            \nabla^i_z 
            \left(
                c_{\Delta}(\cdot, \xi_{\Delta})
                -
                \hat{c}_{\Delta}(\cdot, \xi_{\Delta})
            \right)
        }_{\infty}
        \\\ \nonumber 
        &\leq
        \sum_{z \in \Lambda_n} \sum_{\Delta \nsubseteq \Tilde{\Lambda}_n}\sum_{i=1}^q\sum_{\xi_{\Delta}}
        \left(
            \norm{\nabla^i_z c_{\Delta}(\cdot, \xi_{\Delta})}_{\infty} + \norm{\nabla^i_z \hat{c}_{\Delta}(\cdot, \xi_{\Delta})}_{\infty}
        \right),
    \end{align}
    where we used that, by Proposition~\ref{proposition:oscillation-equation}, it holds that for all $\eta \in \Omega$ 
    \begin{align*}
        \sum_{\Delta \subset \Tilde{\Lambda}_n} \sum_{\xi_{\Delta}} 
            \nabla^i_z 
            \left(
                c_{\Delta}(\cdot, \xi_{\Delta})
                -
                \hat{c}_{\Delta}(\cdot, \xi_{\Delta})
            \right)(\eta)
        =
        -
        \sum_{\Delta \nsubseteq \Tilde{\Lambda}_n}
        \sum_{\xi_{\Delta}} 
            \nabla^i_z 
            \left(
                c_{\Delta}(\cdot, \xi_{\Delta})
                -
                \hat{c}_{\Delta}(\cdot, \xi_{\Delta})
            \right)(\eta). 
    \end{align*}
    By translation invariance of the rates we have 
    \begin{align*}
        \norm{\nabla^i_z c_{\Delta}(\cdot, \xi_{\Delta})}_\infty = \norm{\nabla^i_{0}c_{\Delta-z}(\cdot, \xi_{\Delta - z})}_\infty 
    \end{align*}
    and 
    \begin{align*}
        \norm{\nabla^i_z \hat{c}_{\Delta}(\cdot, \xi_{\Delta})}_\infty = \norm{\nabla^i_{0}\hat{c}_{\Delta-z}(\cdot, \xi_{\Delta - z})}_\infty.
    \end{align*}
    So by a change of variable, $\Xi = \Delta - z$, we get a sum over all $\Delta \Subset \Z^d$ with $\Delta \nsubseteq B_{m-1}:=[-m+1,m-1]^d$ when 
    \begin{align}
        z \in \Lambda_{n, m}:= [-2^n + n + m +1, 2^n - n - m -1]^d, 
    \end{align}
    and for $z \notin \Lambda_{n,m}$ we can just estimate the sum by the sum over all $\Delta \Subset \Z^d$. This gives us 
    \begin{align}
        \abs{\mathbf{III}} 
        \leq 
        \norm{\varphi_n}_{\infty} 
        &\leq 
        \abs{\Lambda_{n, m}}\sum_{\Delta \Subset \Z^d: \ \Delta \nsubseteq B_{m-1}} \sum_{i=1}^q \sum_{\xi_{\Delta}}
            \left(
                \norm{\nabla^i_0 c_{\Delta}(\cdot, \xi_{\Delta})}_{\infty}
                +
                \norm{\nabla^i_0 \hat{c}_{\Delta}(\cdot, \xi_{\Delta})}_{\infty} 
            \right)
        \nonumber \\\ 
        &\quad +\abs{\Lambda_n \setminus \Lambda_{n,m}}\sum_{\Delta \Subset \Z^d} \sum_{i=1}^q \sum_{\xi_{\Delta}}
            \left(
                \norm{\nabla^i_0 c_{\Delta}(\cdot, \xi_{\Delta})}_{\infty}
                +
                \norm{\nabla^i_0 \hat{c}_{\Delta}(\cdot, \xi_{\Delta})}_{\infty} 
            \right).
    \end{align}
    If we fix $m$ and take the density limit, then we obtain 
    \begin{align}
        0 \leq \limsup_{n \to \infty}\frac{1}{\abs{\Lambda_n}}\abs{\mathbf{III}} 
        \leq 
        \sum_{\Delta \Subset \Z^d: \ \Delta \nsubseteq B_{m-1}} \sum_{i=1}^q \sum_{\xi_{\Delta}}
        \left(
            \norm{\nabla^i_0 c_{\Delta}(\cdot, \xi_{\Delta})}_{\infty}
            +
            \norm{\nabla^i_0 \hat{c}_{\Delta}(\cdot, \xi_{\Delta})}_{\infty} 
        \right),
    \end{align}
    because $\abs{\Lambda_n \setminus \Lambda_{n,m}} = o(\abs{\Lambda_n})$. 
    Since this holds for all $m \in \N$ and the right-hand side converges to $0$ as $m$ tends to infinity by our assumptions and Lemma~\ref{lemma:regularity-time-reversal}, we can finally conclude that $\mathbf{III} = o(\abs{\Lambda_n})$. 
\end{proof}

\noindent 
Combining the estimate of the boundary contributions and the rewriting of the bulk contribution allows us to prove our first main result. 
\newline 

\begin{proof}[Proof of Theorem \ref{proposition:relative-entropy-loss-non-positive}]
By Lemma~\ref{lemma:boundary-contributions-boundary-order} we know that for all $\nu \in \calM_1(\Omega)$
\begin{align*}
    g_\mathscr{L}(\nu| \mu) 
    = \limsup_{n \to \infty} \abs{\Lambda_n}^{-1}g_\mathscr{L}^n(\nu| \mu) 
    &\leq 
    \liminf_{n \to \infty}\abs{\Lambda_n}^{-1}\tilde{g}_\mathscr{L}^n(\nu| \mu)
    \\\
    &\leq 
    \limsup_{n \to \infty}\abs{\Lambda_n}^{-1}\tilde{g}_\mathscr{L}^n(\nu| \mu)
    =
    \tilde{g}_\mathscr{L}(\nu|\mu).
\end{align*}
Now Lemma~\ref{lemma:bulk-contribution} tells us that for all $\nu \in \calM_1(\Omega)$ it holds that
\begin{align*}
    \tilde{g}_\mathscr{L}(\nu|\mu) \leq 0.
\end{align*}
By the fundamental theorem of calculus this also implies that 
\begin{align*}
    h(\nu_t | \mu) - h(\nu | \mu) 
    &= 
    \limsup_{n \to \infty}
        \frac{1}{\abs{\Lambda_n}}
        \left(
        h_{\Lambda_n}(\nu_t | \mu) 
        -
        h_{\Lambda_n}(\nu | \mu)
        \right)
        \\\
   &=
    \limsup_{n \to \infty}
        \frac{1}{\abs{\Lambda_n}}
        \int_{0}^t g_\mathscr{L}^n(\nu_s | \mu) \ ds
    \leq 0. 
\end{align*}
This concludes the proof. 
\end{proof}

\subsection{Proof of Theorem \ref{theorem:dgvp}}\label{section:proof-dgvp}

The main work will be to establish the upper-semicontinuity of the (approximating) relative entropy loss. We will do this in two steps. We first define a tentative approximation of the approximating relative entropy loss which satisfies a monotonicity property that allows us to conclude upper-semicontinuity of the limit. As a second step, we then establish that the approximation error vanishes in the density limit. 

To be precise, for $n \in \N$ and $\nu \in \calM_1(\Omega)$ we define  
\begin{align*}
    s_n(\nu|\mu) := \sum_{\eta \in \Omega_{\Lambda_n}}\sum_{\Delta \subset \Tilde{\Lambda}_n}\sum_{\xi_{\Delta} \neq \eta_{\Delta}}f(\nu, \Lambda_n, \eta, \xi_{\Delta})\hat{c}^{B_{n-1}(x(\Delta))}_{\Delta}(\xi_{\Delta}\eta_{\Delta^c},\eta_{\Delta}),
\end{align*}
where we use the monotone truncation
\begin{align*}
    \hat{c}^{\Lambda}_{\Delta}(\eta, \xi_{\Delta}) 
    := 
    \inf_{\omega \in \Omega}c_{\Delta}(\eta_{\Lambda} \omega_{\Lambda^c}, \xi_{\Delta}), 
    \quad 
    B_{n}(x)
    := 
    \set{y \in \Z^d: \abs{x-y}\leq n}.
\end{align*}
For the finitely many distinct basic shapes $\Delta_1, \dots, \Delta_M$ with $0 \in \Delta_i$ the centers of the balls around which we truncate are chosen to be $x(\Delta_i)=0$. For translations $\Delta = \Delta_i + z$, $z \in \Z^d$, of these basic shapes we take $x(\Delta)=z$. 
In the following, we will often just write $x$ instead of $x(\Delta)$ whenever it is clear from the context.
Moreover, we also approximate the function $F$ from the previous section by 
\begin{align*}
    &f(\nu, \Lambda_n, \eta_{\Lambda_n}, \xi_{\Delta}\eta_{\Lambda_n \setminus \Delta}) 
    \\\
    &:= 
    \begin{cases}
        F_0\left(
            \frac{1}{\nu(\xi_{\Delta}\eta_{\Lambda_n \setminus \Delta})}\int\mathbf{1}_{\eta_{\Lambda_n}}(\omega)\frac{\gamma_{\Delta}(\eta_{\Delta}|\omega_{\Delta^c})}{\gamma_{\Delta}(\xi_{\Delta}|\omega_{\Delta^c})}\nu(d\omega)
        \right)
        \nu(\xi_{\Delta}\eta_{\Lambda_n \setminus \Delta}), 
        \text{ if $\nu(\xi_{\Delta}\eta_{\Lambda_n \setminus \Delta})>0$,} 
        \\\
        -\infty,  
        \hspace{5.6cm} \text{ if $\nu(\xi_{\Delta}\eta_{\Lambda_n \setminus \Delta})= 0$  and $\nu(\eta_{\Lambda_n})>0$},
        \\\
        0,  
        \hspace{6.1cm}\text{ if $\nu(\xi_{\Delta}\eta_{\Lambda_n \setminus \Delta})= \nu(\eta_{\Lambda_n})=0$}. 
    \end{cases}
\end{align*}
\noindent 
We show the existence of the density limit of $(s_n(\cdot |\mu))_{n \in \N}$ on $\calM^{inv}_1(\Omega)$ in two steps. First we show that it satisfies a growth property that is reminiscent of subadditivity and afterwards we combine this with a multiplicative volume correction to show the convergence via a monotonicity argument.

\begin{lemma}\label{lemma:growth-bound-approximation}
    Assume that the rates satisfy conditions $\mathbf{(R2)}$ and $\mathbf{(R4)}$ and that $\mu \in \mathscr{G}(\gamma)$ is translation-invariant, where the specification $\gamma$ satisfies condition $\mathbf{(S4)}$. 
    Let $m \in \N$ be such that for all $\Delta \Subset \Z^d$ with $0 \in \Delta$ and $c_\Delta > 0$ we have $\Delta \subset B_{m-1}(0)$.
    Then, for all $n \geq m$ and all $\nu \in \calM^{inv}_1(\Omega)$, it holds that 
    \begin{align}
        s_n(\nu|\mu) \leq 2^d s_{n-1}(\nu|\mu).
    \end{align}
\end{lemma}

\begin{proof}
    The main argument in the proof is to upper bound $s_n(\nu|\mu)$ by $s_{n-1}(\nu|\mu)$ using Jensen's inequality and the concavity of $F_0$. 
    To make this precise, consider $2^d$ disjoined and congruent subcubes $\Gamma_{n,k}$ of $\Lambda_n$ with total side-length $2^n -1$ as well as $2^d$ disjoined and congruent subcubes $\Tilde{\Gamma}_{n,k}$ of $\Tilde{\Lambda}_n$ with total side-length $2^n -n -1$. Let the subcubes be centered such that $\Tilde{\Gamma}_{n,k} \subset \Gamma_{n,k}$ for each $k$. Note that $\cup_{k=1}^{2^d}\Tilde{\Gamma}_{n,k} \subsetneq \Tilde{\Lambda}_n$. 
    Now we can estimate 
    \begin{align*}
        s_n(\nu |\mu) 
        &= 
        \sum_{\eta_{\Lambda_n}}\sum_{\Delta \subset \Tilde{\Lambda}_n}\sum_{\xi_{\Delta} \neq \eta_{\Delta}}
            f(\nu, \Lambda_n, \eta, \xi_{\Delta})
            \hat{c}^{B_{n-1}(x(\Delta))}_{\Delta}(\xi_{\Delta}\eta_{\Delta^c} \eta_{\Delta})
        \nonumber \\\
        &\leq 
        \sum_{k=1}^{2^d}\sum_{\eta_{\Lambda_n}}\sum_{\Delta \subset \Tilde{\Lambda}_n: \ x(\Delta) \in \Tilde{\Gamma}_{n,k}}\sum_{\xi_{\Delta} \neq \eta_{\Delta}}
            f(\nu, \Lambda_n, \eta, \xi_{\Delta})
            \hat{c}^{B_{n-1}(x(\Delta))}_{\Delta}(\xi_{\Delta}\eta_{\Delta^c}, \eta_{\Delta}),
    \end{align*}
    where the inequality is due to the fact that the subcubes do not cover all of $\Tilde{\Lambda}_n$ and $f$ is non-positive. Now for each $k=1,\dots, 2^d$ we can split the summation over $\eta_{\Lambda_n}$ into two summations, one over $\eta_{\Lambda_n \setminus \Gamma_{n,k}}$ and one over $\omega_{\Gamma_{n,k}}$. This gives us
    \begin{align}\label{proof:growth-approximation-1}
        \nonumber 
        &\sum_{k=1}^{2^d}\sum_{\eta_{\Lambda_n \setminus \Gamma_{n,k}}} \sum_{\omega_{\Gamma_{n,k}}}\sum_{\Delta \subset \Tilde{\Lambda}_n: \ x(\Delta) \in \Tilde{\Gamma}_{n,k}}\sum_{\xi_{\Delta} \neq \eta_{\Delta}}
            f(\nu, n, \omega_{\Gamma_{n,k}}\eta_{\Lambda_n \setminus \Gamma_{n,k}}, \xi_{\Delta})
            \nonumber 
            \\\
            &\hspace{5.5cm}
            \times \hat{c}^{B_{n-1}(x(\Delta))}_{\Delta}(\xi_{\Delta}(\omega\eta)_{\Delta^c}, (\omega\eta)_{\Delta})
         \\\
        = 
        &\sum_{k=1}^{2^d}\sum_{\eta_{\Gamma_{n,k}}}\sum_{\Delta \subset \Tilde{\Lambda}_n: \ x(\Delta) \in \Tilde{\Gamma}_{n,k}}\sum_{\xi_{\Delta} \neq \eta_{\Delta}}
        \hat{c}^{B_{n-1}(x(\Delta))}_{\Delta}(\xi_{\Delta}\eta_{\Delta^c}, \eta_{\Delta})
        \sum_{\eta_{\Lambda_n \setminus \Gamma_{n,k}}} 
        f(\nu, n, \eta_{\Lambda_n}, \xi_{\Delta}),
        \nonumber
    \end{align}
    where we were able to pull the rates out of the summation, because for $n$ sufficiently large we have $\Delta \subset \Gamma_{n,k}$ if $x(\Delta) \in \Tilde{\Gamma}_{n,k}$ by the assumption that there are only finitely many distinct basic shapes $\Delta$ on which we can perform updates and by construction of the subcubes we also have $B_{n-1}(x(\Delta)) \subset \Gamma_{n,k}$ for $x(\Delta) \in \Tilde{\Gamma}_{n,k}$. 
    Since $F_0$ is a concave function, we can use the definition of $f$ and Jensen's inequality to get
    \begin{align*}
        \sum_{\eta_{\Lambda_n \setminus \Gamma_{n,k}}} 
        f(\nu, n, \eta_{\Lambda_n}, \xi_{\Delta})
        &= 
        \nu(\xi_{\Delta}\eta_{\Gamma_{n,k}\setminus \Delta})
        \sum_{\eta_{\Lambda_n \setminus \Gamma_{n,k}}}
        \frac{\nu(\xi_{\Delta}\eta_{\Lambda_n \setminus \Delta})}
        {\nu(\xi_{\Delta}\eta_{\Gamma_{n,k}\setminus \Delta})}
        \\\
        &\hspace{1cm}
        \times 
        F_0\left(
            \frac{1}{\nu(\xi_{\Delta}\eta_{\Lambda_n \setminus \Delta})}\int\mathbf{1}_{\eta_{\Lambda_n}}(\omega)\frac{\gamma_{\Delta}(\eta_{\Delta}|\omega_{\Delta^c})}{\gamma_{\Delta}(\xi_{\Delta}|\omega_{\Delta^c})}\nu(d\omega)
        \right)
        \nonumber \\\
        \leq 
        \nu(\xi_{\Delta}\eta_{\Gamma_{n,k}\setminus \Delta})
        &F_0\left(
            \sum_{\eta_{\Lambda_n \setminus \Gamma_{n,k}}}
                \frac{1}
                {\nu(\xi_{\Delta}\eta_{\Gamma_{n,k}\setminus \Delta})}
                \int\mathbf{1}_{\eta_{\Lambda_n}}(\omega)\frac{\gamma_{\Delta}(\eta_{\Delta}|\omega_{\Delta^c})}{\gamma_{\Delta}(\xi_{\Delta}|\omega_{\Delta^c})}\nu(d\omega)
            \right)
         \\\
        =
        \nu(\xi_{\Delta}\eta_{\Gamma_{n,k}\setminus \Delta})
        &F_0\left(
                \frac{1}
                {\nu(\xi_{\Delta}\eta_{\Gamma_{n,k}\setminus \Delta})}
                \int\mathbf{1}_{\eta_{\Lambda_n}}(\omega)\frac{\gamma_{\Delta}(\eta_{\Delta}|\omega_{\Delta^c})}{\gamma_{\Delta}(\xi_{\Delta}|\omega_{\Delta^c})}\nu(d\omega)
            \right). 
        \nonumber
    \end{align*}
    Plugging this back into~\eqref{proof:growth-approximation-1} and using the translation invariance of the rates and $\nu$, implies that
    \begin{align*}
        s_n(\nu|\mu) 
        &\leq 
        \sum_{k=1}^{2^d}\sum_{\eta_{\Gamma_{n,k}}}\sum_{\Delta \subset \Tilde{\Lambda}_n : \ x(\Delta) \in \Tilde{\Gamma}_{n,k}}\sum_{\xi_{\Delta} \neq \eta_{\Delta}}
        \hat{c}^{B_{n-1}(x(\Delta))}_{\Delta}(\xi_{\Delta}\eta_{\Delta^c}, \eta_{\Delta})
        \nu(B(\Gamma_{n,k}, \xi_{\Delta}\eta_{\Gamma_{n,k}\setminus \Delta}))
        \nonumber \\\
        &\hspace{4cm}
        \times 
        F_0\left(
                \frac{1}
                {\nu(\xi_{\Delta}\eta_{\Gamma_{n,k}\setminus \Delta})}
                \int\mathbf{1}_{\eta_{\Gamma_{n,k}}}(\omega)\frac{\gamma_{\Delta}(\eta_{\Delta}|\omega_{\Delta^c})}{\gamma_{\Delta}(\xi_{\Delta}|\omega_{\Delta^c})}\nu(d\omega)
            \right) 
        \\\
        &\leq 
        2^d s_{n-1}(\nu|\mu).
        \nonumber 
    \end{align*}
    Note that for the last inequality we also used that truncating the rates over a smaller volume is non-increasing and that the function $f$ is non-positive by definition. This allowed us to drop some terms from the summation without any harm. 
\end{proof}

\noindent 
With this growth property it is now easy to conclude the convergence of $(\abs{\Lambda_n}^{-1}s_n(\nu| \mu))_{n \in \N}$ and the upper-semicontinuity of the limit. 
\begin{lemma}\label{lemma:convergence-of-approximation}
    Assume that the rates satisfy conditions $\mathbf{(R2)}$ and $\mathbf{(R4)}$ and that $\mu \in \mathscr{G}(\gamma)$ is translation-invariant, where the specification $\gamma$ satisfies condition $\mathbf{(S4)}$.
    Then, the following limit exists for all $\nu \in \calM_1^{inv}(\Omega)$ and defines an upper-semicontinuous function on $\calM_1^{inv}$: 
    \begin{align}
        s(\nu|\mu):=\lim_{n \to \infty}\frac{1}{\abs{\Lambda_n}}s_n(\nu|\mu). 
    \end{align}
\end{lemma}
\begin{proof}
    The previous estimate on the growth of $s_n(\nu| \mu)$ looks very similar to classical subadditivity of $(s_n(\nu|\mu))_{n \in \N}$, which would be sufficient to conclude existence and upper-semicontinuity of the limit via a $d$-dimensional generalisation of Fekete's Lemma, see~\cite[Lemma 15.11]{georgii_gibbs_2011}. However, we cannot apply this result directly, but have to account for the volume growth of $\Lambda_n$ via a multiplicative correction that goes to $1$ as $n$ tends to infinity. More precisely,
    for $n \in \N$ we define the volume correction
    \begin{align*}
        G_n := \prod_{k=n}^{\infty}\frac{(2^{k+2}-2)^d}{(2^{k+2}-1)^d}.
    \end{align*}
    Then by the Lemma \ref{lemma:growth-bound-approximation} we have
    \begin{align*}
        \frac{G_{n+1}}{\abs{\Lambda_{n+1}}}s_{n+1}(\nu|\mu) \leq \frac{G_{n}}{\abs{\Lambda_{n}}}s_{n}(\nu|\mu). 
    \end{align*} 
    Therefore, the limit 
    \begin{align*}
        \lim_{n\to \infty}\frac{G_{n}}{\abs{\Lambda_{n}}}s_{n}(\nu|\mu)
    \end{align*}
    exists by monotonicity and is upper-semicontinous as the limit of a non-increasing sequence of upper-semicontinuous functions. 
    Since $G_n$ converges to $1$ as $n$ tends to infinity, this implies that the following limit also exists and is equal to the above:
    \begin{align*}
        \lim_{n\to \infty}\frac{1}{\abs{\Lambda_{n}}}s_{n}(\nu|\mu) = \lim_{n\to \infty}\frac{G_{n}}{\abs{\Lambda_{n}}}s_{n}(\nu|\mu),
    \end{align*}
    as desired.
\end{proof}

As a second step we now show that the approximation error we make vanishes in the density limit. This shows in particular that the approximating relative entropy loss functional is upper-semicontinuous on $\calM_1^{inv}(\Omega)$.
More precisely, because of Lemma~\ref{lemma:bulk-contribution} we need to show that $s_n(\nu|\mu)$ is really an approximation to
\begin{align}
    S_n(\nu|\mu) 
    := 
    \sum_{\eta_{\Lambda_n}}\sum_{\Delta \subset \Tilde{\Lambda}_n}\sum_{\xi_{\Delta} \neq \eta_{\Delta}} 
        F(\nu, n, \eta_{\Lambda_n}, \xi_{\Delta})c^{(n)}_{\Delta}(\eta,\xi_{\Delta})
        \frac{\mu(\eta_{\Lambda_n})}{\mu(\xi_{\Delta}\eta_{\Lambda_n \setminus \Delta})}
    = \Tilde{g}_{\mathscr{L}}^n(\nu | \mu) + o(\abs{\Lambda_n}). 
\end{align}
\noindent 
To show that the approximation error we make by replacing $S_n$ with $s_n$ is of boundary order, and therefore negligible when taking the density limit, we will also reuse Lemma~\ref{lemma:uniform-convergence-conditional-probabilities}. 

\begin{lemma}\label{lemma:approximation-is-approximation}
    Assume that the rates $(c_{\Delta}(\cdot, \xi_{\Delta}))_{\Delta \Subset \Z^d, \xi_{\Delta}\in \Omega_\Delta}$ satisfy the conditions $\mathbf{(R1)}$, $\mathbf{(R2)}$, $\mathbf{(R4)}$, and $\mathbf{(R5)}$. 
    Moreover, assume that $\mu \in \mathscr{G}(\gamma)$ is time-stationary and translation-invariant, where the specification $\gamma$ satisfies conditions $\mathbf{(S1)}$, $\mathbf{(S2)}$, and $\mathbf{(S4)}$. 
    Then, the density limit of $(S_n(\nu|\mu))_{n \in \N}$ exists and is equal to $s(\nu|\mu)$, i.e., 
    \begin{align}
        \Tilde{g}_\mathscr{L}(\nu|\mu)
        =
        \lim_{n \to \infty}\frac{1}{\abs{\Lambda_n}}S_n(\nu|\mu) = s(\nu|\mu). 
    \end{align}
    If the rates additionally satisfy condition $\mathbf{(R3)}$ and the specification also satisfies $\mathbf{(S3)}$, then the density limit of $(S_n(\nu|\mu))_{n \in \N}$ agrees with the density limit of $(\Tilde{g}_\mathscr{L}(\nu|\mu))_{n \in \N}$.
    In particular, the approximating entropy loss per site $\Tilde{g}_\mathscr{L}(\cdot|\mu)$ is then an upper-semicontinuous functional on $\calM_1^{\text{inv}}(\Omega)$. 
\end{lemma} 

\begin{proof}
    Let $n \in \N$. 
    If there was an $\eta \in \Omega_{\Lambda_n}$ such that $\nu(\xi_{\Delta}\eta_{\Lambda_n \setminus \Delta})=0$ and $\nu(\eta_{\Lambda_n})>0$, 
    then we have $S_m(\nu|\mu) = s_m(\nu|\mu)$ in the sense of $-\infty = -\infty$ for all $m \geq n$. 
    If there was an $\eta \in \Omega_{\Lambda_n}$ such that $\nu(\xi_{\Delta}\eta_{\Lambda_n \setminus \Delta})=0$ and $\nu(\eta_{\Lambda_n})=0$, 
    then the corresponding summands would not contribute to either summation. 
    Therefore, we can assume without loss of generality that $\nu(\eta_{\Lambda_n})>0$ for all $\eta \in \Omega_{\Lambda_n}$. 
    In this case, both $s_n(\nu | \mu)$ and $S_n(\nu | \mu)$ are finite, and we can write 
    \begin{align}
        S_n(\nu | \mu)& - s_n(\nu | \mu) 
        = 
        \sum_{\eta_{\Lambda_n}}\sum_{\Delta \subset \Tilde{\Lambda}_n}\sum_{\xi_{\Delta} \neq \eta_{\Delta}}
            \Bigg[
                F(\nu,n,\eta_{\Lambda_n}, \xi_{\Delta})c^{(n)}_{\Delta}(\eta, \xi_{\Delta})
                \frac{\mu(\eta_{\Lambda_n})}{\mu(\xi_{\Delta}\eta_{\Lambda_n \setminus \Delta})}
                \nonumber 
                \\\
                &\hspace{2.5cm}- 
                f(\nu, n, \eta{\Lambda_n}, \xi_{\Delta})
                \frac{c^{(n)}_{\Delta}({\Lambda_n}, \xi_{\Delta}) \nu(\eta_{\Lambda_n})}
                {
                \int
                    \mathbf{1}_{\eta_{\Lambda_n}}(\omega)
                    \frac{\gamma_{\Delta}(\xi_{\Delta}|\omega_{\Delta^c})}
                    {\gamma_{\Delta}(\eta_{\Delta}|\omega_{\Delta^c})}
                \nu(d\omega)
                }
            \Bigg]
            \nonumber 
        \\\ 
        &-
        \sum_{\eta_{\Lambda_n}}\sum_{\Delta \subset \Tilde{\Lambda}_n}\sum_{\xi_{\Delta} \neq \eta_{\Delta}}
            f(\nu,n, \eta_{\Lambda_n}, \xi_{\Delta})
            \left[
            \hat{c}^{B_{n-1}(x)}_{\Delta}(\xi_{\Delta}\eta_{\Delta^c}, \eta_{\Delta})
            -
            \tfrac{c^{(n)}_{\Delta}(\eta_{\Lambda_n}, \xi_{\Delta}) \nu(\eta_{\Lambda_n})}
                {
                \int
                    \mathbf{1}_{\eta_{\Lambda_n}}(\omega)
                    \frac{\gamma_{\Delta}(\xi_{\Delta}|\omega_{\Delta^c})}
                    {\gamma_{\Delta}(\eta_{\Delta}|\omega_{\Delta^c})}
                \nu(d\omega)
                }
            \right]
            \nonumber 
            \\\
            &=: \mathbf{I} + \mathbf{II}. 
        \nonumber 
    \end{align}
    We will now proceed by estimating these two terms separately. We start with $\mathbf{I}$ and use the definitions of the functions $F$ and $f$ to obtain
    \begin{align}
        \abs{\mathbf{I}}
        \leq 
        \sum_{\eta_{\Lambda_n}}\sum_{\Delta \subset \Tilde{\Lambda}_n}\sum_{\xi_{\Delta} \neq \eta_{\Delta}}
            &\nu(\eta_{\Lambda_n})c^{(n)}_{\Delta}(\eta_{\Lambda_n}, \xi_{\Delta})
            \nonumber 
            \\\
            &\times \abs{
                \log\left( 
                    \frac{\mu(\eta_{\Lambda_n})}{\mu(\xi_{\Delta}\eta_{\Lambda_n \setminus \Delta})}
                    \frac{1}{\nu(\eta_{\Lambda_n})}
                    \int
                        \mathbf{1}_{\eta_{\Lambda_n}}(\omega)
                        \frac{\gamma_{\Delta}(\xi_{\Delta}|\omega_{\Delta^c})}{\gamma_{\Delta}(\eta_{\Delta}|\omega_{\Delta^c})}
                    \nu(d\omega)
                \right)
            }.
    \end{align}
    The terms inside the logarithm converge to $1$ uniformly in $\eta$ by Lemma ~\ref{lemma:uniform-convergence-conditional-probabilities}, where we again use that the convergence is moreover uniform for all summands by assumptions $\mathscr{(R2)}$ and $\mathbf{(R4)}$.   
    In combination with the combinatorial estimate from Lemma \ref{combinatorial-lemma} this implies that $\abs{\mathbf{I}} = o(\abs{\Lambda_n})$.

    To deal with $\mathbf{II}$, we first note that by Lemma \ref{differentiation-lemma} we have
    \begin{align}\label{estimate-bracket-term}
        \lim_{n \to \infty} 
            \frac{1}{\nu(\eta_{\Lambda_n})} 
            \int
                \mathbf{1}_{\eta_{\Lambda_n}}(\omega)
                \frac{\gamma_{\Delta}(\xi_{\Delta}|\omega_{\Delta^c})}{\gamma_{\Delta}(\eta_{\Delta}|\omega_{\Delta^c})}
            \nu(d\omega)
            = \frac{\gamma_{\Delta}(\xi_{\Delta}|\eta_{\Delta^c})}{\gamma_{\Delta}(\eta_{\Delta}|\eta_{\Delta^c})},
    \end{align}
    and the convergence holds uniform in $\eta \in \Omega$ by quasilocality of the specification $\gamma$. Now we would like to conclude that $\mathbf{II} = s_n(\nu | \mu) o(1)$, which would yield $\mathbf{II}=o(\abs{\Lambda_n})$ by convergence of $(s_n(\nu |\mu))_{n \in \N}$. But for this, we have to make sure that we are actually allowed to divide by $\hat{c}_\Delta^{B_{n-1}(x)}(\xi_\Delta \eta_{\Delta^c}, \eta_\Delta)$. As we will see in Lemma~\ref{lemma:dichotomy-helper} after the end of this proof, there exists $N \in \N$, uniform in $\Delta$ and $\xi_{\Delta}$ and $\eta \in \Omega$, such that for each fixed $\eta \in \Omega$ we either have 
    \begin{align*}
        \hat{c}_\Delta^{B_{n-1}(x)}(\xi_\Delta \eta_{\Delta^c}, \eta_\Delta) = 0
    \end{align*}
    for all $n \in \N$ or 
    \begin{align*}
        \hat{c}_\Delta^{B_{n-1}(x)}(\xi_\Delta \eta_{\Delta^c}, \eta_\Delta) > 0
    \end{align*}
    for all $n \geq N$. In the former case, the terms corresponding to such $\eta$ do not enter the summation for any $n \in \N$ and in the latter case we are allowed to perform the division for sufficiently large $n$. Therefore, we can assume without loss of generality that 
    $
        \hat{c}_\Delta^{B_{n-1}(x)}(\xi_\Delta \eta_{\Delta^c}, \eta_\Delta) > 0
    $
    for all $\xi_{\Delta}$ and $\eta$. This allows us to conclude
    \begin{align*}
        \mathbf{II} 
        &= 
        \sum_{\eta_{\Lambda_n}}\sum_{\Delta \subset \Tilde{\Lambda}_n}\sum_{\xi_{\Delta} \neq \eta_{\Delta}}
        f(\nu,n, \eta_{\Lambda_n}, \xi_{\Delta})
        \hat{c}^{B_{n-1}(x)}_{\Delta}(\xi_{\Delta}\eta_{\Delta^c}, \eta_{\Delta})
        \nonumber 
        \\\
        &\hspace{2.5cm}
        \times 
        \left[
        1
        -
        \frac{1}{\hat{c}^{B_{n-1}(x)}_{\Delta}(\xi_{\Delta}\eta_{\Delta^c}, \eta_{\Delta})}
        \frac{c^{(n)}_{\Delta}(\eta_{\Lambda_n}, \xi_{\Delta}) \nu(\eta_{\Lambda_n})}
            {
            \int
                \mathbf{1}_{\eta_{\Lambda_n}}(\omega)
                \frac{\gamma_{\Delta}(\xi_{\Delta}|\omega_{\Delta^c})}
                {\gamma_{\Delta}(\eta_{\Delta}|\omega_{\Delta^c})}
            \nu(d\omega)
            }
        \right]
         \\\ 
        &= s_n(\nu|\mu) \ o(1)= o(\abs{\Lambda_n}), 
        \nonumber
    \end{align*}
    as desired.
\end{proof}

\noindent 
At the end of the previous proof, we used the following rather technical lemma to make sure that we are not performing a division by zero. 

\begin{lemma}\label{lemma:dichotomy-helper}
    Assume that the rates $(c_{\Delta}(\cdot, \xi_{\Delta}))_{\Delta \Subset \Z^d, \xi_{\Delta}}$ satisfy the conditions $\mathbf{(R1)}$, $\mathbf{(R2)}$, $\mathbf{(R4)}$, and $\mathbf{(R5)}$.
    Then, there exists $N \in \N$ such that for all $\eta \in \Omega$ we have the following dichotomy for the truncated rates: 
    \begin{align*}
        \forall \Delta \Subset \Z^d \ \forall \xi_{\Delta} : 
        \quad 
        \text{Either }
        c^{B_{N-1}(x(\Delta))}_{\Delta}(\eta_{\Lambda_N}, \xi_{\Delta}) >0 
        \text{ or }
        \forall n \in \N: \ c^{B_{n-1}(x(\Delta))}_{\Delta}(\eta_{\Lambda_n}, \xi_{\Delta}) =0.  
    \end{align*}
    Moreover, if the rates $(c_{\Delta}(\cdot, \xi_{\Delta}))_{\Delta \Subset \Z^d, \xi_{\Delta}}$ satisfy the above conditions and the specification $\gamma$ satisfies $\mathbf{(S1)}$, $\mathbf{(S2)}$,  and $\mathbf{(S4)}$,
    then the rates of the time-reversal $(\hat{c}_{\Delta}(\cdot, \xi_{\Delta}))_{\Delta \Subset \Z^d, \xi_{\Delta}}$ also satisfy $\mathbf{(R1)}$, $\mathbf{(R2)}$, $\mathbf{(R4)}$, and $\mathbf{(R5)}$. 
\end{lemma}

\begin{proof}
    By translation invariance of the rates, we only have to worry about those $\Delta \Subset \Z^d$ with $0 \in \Delta$ and $c_{\Delta}>0$. 
    Let $\kappa:=  \inf_{\Delta \Subset \Z^d, \xi_{\Delta}, \eta: \ c_{\Delta}(\eta, \xi_{\Delta})>0}c_{\Delta}(\eta, \xi_{\Delta}) > 0$ be the minimal transition rate.     
    Since the local state space is finite, the continuity of the rates implies that they are also quasilocal. 
    Since we also assumed that there are only finitely many types of transitions, there exists $N \in \N$, uniform in $\Delta \Subset \Z^d$, such that if two configurations agree on $B_{N-1}(x(\Delta))$,
    then we have 
    \begin{align*}
        \abs{c_{\Delta}(\omega_1, \xi_{\Delta})-c_{\Delta}(\omega_2, \xi_{\Delta})} < \frac{\kappa}{2}. 
    \end{align*}
    In particular, this implies that if $\eta \in \set{\omega\colon c_{\Delta}(\omega, \xi_{\Delta})>0}$, then 
    \begin{align*}
        c^{B_{N-1}(x(\Delta))}_{\Delta}(\eta_{\Lambda_N}, \xi_{\Delta}) >0,
    \end{align*}
    and if $\eta \notin \set{\omega\colon c_{\Delta}(\omega, \xi_{\Delta})>0}$, then 
    \begin{align*}
        \forall n \in \N: \ c^{B_{n-1}(x(\Delta))}_{\Delta}(\eta_{\Lambda_n}, \xi_{\Delta}) =0.
    \end{align*}
    To see that the rates $\hat{c}$ of the time-reversal also satisfy the conditions we just combine the corresponding assumptions on the rates $c$ and on the specification $\gamma$. 
\end{proof}

As a final ingredient for the proof of Theorem \ref{theorem:dgvp} we need to show that, if a measure has vanishing approximating entropy loss with respect to $\mu$, then it is itself a Gibbs measure with respect to the specification $\gamma$. 
Note that we use the irreducibility assumption $\mathbf{(R6)}$ for the first time here. 

\begin{lemma}\label{lemma:zero-loss-gibbs}
    Assume that the rates satisfy conditions $\mathbf{(R1)}-\mathbf{(R6)}$ and that $\mu \in \mathscr{G}(\gamma)$ is time-stationary for the dynamics, where the specification $\gamma$ satisfies conditions $\mathbf{(S1)}-\mathbf{(S4)}$. 
    Let $\nu \in \calM^{inv}_1(\Omega)$. If $\tilde{g}_{\mathscr{L}}(\nu|\mu)=0$, then $\nu \in \mathscr{G}(\gamma)$. 
\end{lemma}
\noindent 
The proof is very similar to the middle part of the proof of \cite[Theorem 2.12]{jahnel_attractor_2019}, but treats a more general situation, since we do not need to use any reversibility assumption. This is because we already eliminated all the dangerous terms in the proof of Lemma \ref{lemma:bulk-contribution} by using Proposition \ref{proposition:oscillation-equation}.

\begin{proof}
    Since the convergence is monotone, the assumption that $\tilde{g}_{\mathscr{L}}(\nu \mid \mu)=0$ implies that we must already have 
    \begin{align*}
         \quad s_n(\nu|\mu) = 0
    \end{align*}
    for all $n \in \N$. 
    By definition of $s_n(\cdot)$, $f$ and $F_0$, this implies that all the terms of the sum have to vanish. Hence, for all $n \in \N$, $\eta_{\Lambda_n} \in \Omega_{\Lambda_n}$, and $j=1,\dots, q$, 
    such that $\hat{c}^{B_{n-1}(x)}_{\Delta}(\xi_{\Delta}\eta_{\Delta^c}, \eta_{\Delta}) > 0$ we either have
    \begin{align}\label{general-key-equality}
        \frac{1}{\nu(\xi_{\Delta}\eta_{\Lambda_n \setminus \Delta})} 
        \int
            \mathbf{1}_{\eta_{\Lambda_n}}(\omega)
            \frac{\gamma_{\Delta}(\xi_{\Delta}|\omega_{\Delta^c})}{\gamma_{\Delta}(\eta_{\Delta}|\omega_{\Delta^c})}
        \nu(d\omega) 
        = 
        1, 
    \end{align}
    or 
    \begin{align*}
        \nu(\eta_{\Lambda_n}) = \nu(\xi_{\Delta}\eta_{\Lambda_n \setminus \Delta}) = 0. 
    \end{align*}
    But if there was $\eta_{\Lambda_n} \in \Omega_{\Lambda_n}$ such that $\nu(\eta_{\Lambda_n}) = 0$, then by irreducibility we would necessarily have
    \begin{align*}
        \forall \psi \in \Omega_{\Tilde{\Lambda}_m\setminus \Lambda_n} \ \forall \Delta \subset \Lambda_n \ \forall \xi_{\Delta}:  \quad \nu(\xi_{\Delta}\eta_{\Lambda_n\setminus \Delta}\psi_{\Tilde{\Lambda}_m \setminus \Lambda_n}) = 0,
    \end{align*}
    for all large enough $m$ such that $\Lambda_n \subset \Tilde{\Lambda}_m$. 
    Since this holds for all boundary conditions $\psi$, all $\Delta \subset \Lambda_n$ and  $\xi_{\Delta}$, we can deduce that
    \begin{align*}
        \forall \eta_{\Lambda_n} \in \Omega_{\Lambda_n}: \quad \nu(\eta_{\Lambda_n})=0, 
    \end{align*}
    which cannot be true. Therefore, we must have 
    \begin{align*}
        \forall n \in \N: \ \forall \eta_{\Lambda_n} \in \Omega_{\Lambda_n}: \quad \nu(\eta_{\Lambda_n}) > 0, 
    \end{align*}
    and in particular~\eqref{general-key-equality} holds if $\hat{c}^{B_{n-1}(x)}_{\Delta}(\xi_{\Delta}\eta_{\Delta^c}, \eta_{\Delta}) > 0$.  
    In this case, we can use martingale convergence and the differentiation lemma to see that by irreducibility
    \begin{align*}
        \nu-\text{a.a. } \eta \in \Omega \ \forall \Delta \Subset \Z^d \ \forall \xi_{\Delta}: \quad 
        \frac{\nu(\eta_{\Delta}|\eta_{\Delta^c})}{\nu(\xi_{\Delta}|\eta_{\Delta^c})}
        =
        \lim_{n \to \infty}\frac{\nu(\eta_{\Lambda_n})}{\nu(\xi_{\Delta}\eta_{\Lambda_n \setminus \Delta})} = \frac{\gamma_{\Delta}(\xi_{\Delta}|\eta_{\Delta^c})}{\gamma_{\Delta}(\eta_{\Delta}|\eta_{\Delta^c})}.
    \end{align*}
    Via the irreducibility assumption and the fact, that if two strictly positive probability vectors $a = (a_1, a_2, ..., a_n)$ and $b = (b_1, ...,b_n)$ satisfy $a_i/a_j = b_i/b_j$ for all $i,j=1,...,n$, then we necessarily have $a=b$, we get 
    \begin{align*}
        \forall \Delta \Subset \Z^d \ \forall \xi_{\Delta}:  
        \quad 
        \gamma_{\Delta}(\xi_{\Delta} \lvert \eta_{\Delta^c})
        = 
        \nu(\xi_{\Delta} \lvert \eta_{\Delta^c}) 
        \quad 
        \text{ for $\nu$-almost all }\eta \in \Omega,
    \end{align*}
    which implies that $\nu \in \mathscr{G}(\gamma)$. 
\end{proof}

\noindent 
Now we have all the ingredients for proving the dynamical Gibbs variational principle for non-reversible interacting particle systems. 
\newline 

\begin{proof}[Proof of Theorem \ref{theorem:dgvp}]
\textit{Ad i. and ii.:} The existence of the limit and its upper-\linebreak semicontinuity follow from Lemma~\ref{lemma:approximation-is-approximation}. 
\newline 
\textit{Ad iii.:} That the approximating relative entropy loss is non-positive was shown in Proposition~\ref{proposition:relative-entropy-loss-non-positive}. 
\newline 
\textit{Ad iv.:} This is exactly what we showed in Lemma \ref{lemma:zero-loss-gibbs}. 
\end{proof}

\subsection{Proof of Theorem \ref{theorem:attractor-property}}\label{section:proof-attractor-property}

We begin with the following technical result. 

\begin{proposition}\label{proposition:attractor-property-helper}
    Assume that the rates $(c_{\Delta}(\cdot, \xi_{\Delta}))_{\Delta \Subset \Z^d, \xi_\Delta\in \Omega_\Delta}$ satisfy conditions $\mathbf{(R1)}-\mathbf{(R5)}$. Moreover, assume that there exists a translation-invariant time-stationary measure $\mu\in\calG(\gamma)$ where the specification $\gamma$ satisfies $\mathbf{(S1)}-\mathbf{(S3)}$. 
    Then, for $\nu \in \calM_1^{inv}(\Omega)$ and $t>0$ we have 
    \begin{align*}
        h(\nu_t \lvert \mu) - h(\nu \lvert \mu) \leq \int_0^t \Tilde{g}_{\mathscr{L}}(\nu_s \lvert \mu)ds. 
    \end{align*}
    Moreover, if $\nu \notin \mathscr{G}(\gamma)$, then there exists a weakly open set $G_{\nu} \subset \calM_1^{inv}(\Omega)$  containing $\nu$, and $\delta, \varepsilon>0$ such that we have 
    \begin{align*}
        \forall \rho  \in G_{\nu} \ \forall 0 \leq s \leq \varepsilon: \quad h(\rho_s| \mu) - h(\rho | \mu) \leq -\delta s . 
    \end{align*}
\end{proposition}

\begin{proof}
    For $n \in \N$ we have seen in Lemma \ref{lemma:boundary-contributions-boundary-order}  that 
    \begin{align*}
        h_{\Lambda_n}(\nu_t\lvert \mu ) - h_{\Lambda_n}(\nu\lvert \mu) \leq \int_0^t \tilde{g}^n_{\mathscr{L}}(\nu_s \lvert \mu)ds + t\cdot o(\abs{\Lambda_n}).
    \end{align*}
    By taking the density limit and using the monotone convergence theorem we obtain
    \begin{align*}
        h(\nu_t \lvert \mu) - h(\nu\lvert \mu) \leq \int_0^t\tilde{g}_{\mathscr{L}}(\nu_s \lvert \mu)ds. 
    \end{align*}
     To prove the second part of the proposition we just use the above formula and the upper-semicontinuity of $\tilde{g}_{\mathscr{L}}(\cdot \lvert \mu)$. 
\end{proof}

\noindent
Now we are ready to state and prove the attractor property for non-reversible interacting particle systems. 
\newline 

\begin{proof}[Proof of Theorem~\ref{theorem:attractor-property}]
    For arbitrary $\nu$ and $n \in \N$ we have by non-nullness of $\mu$ that 
    \begin{align}\label{proof:attractor-property-1}
        \forall \Lambda \Subset \Z^d: \quad -M\abs{\Lambda} \leq h_{\Lambda}(\nu\lvert \mu) \leq M\abs{\Lambda},
    \end{align}
    for some $M>0$ that does not depend on $\nu$, only on $\mu$. 
    Indeed, we can decompose
    \begin{align*}
        h_{\Lambda}(\nu \lvert \mu) 
        = 
        \sum_{\eta_{\Lambda}}\nu(\eta_{\Lambda})\log\left(\nu(\eta_{\Lambda})\right)
        -
        \sum_{\eta_{\Lambda}}\nu(\eta_{\Lambda})\log\left(\mu(\eta_{\Lambda})\right).
    \end{align*}
The first sum is bounded from below by $0$ and from above by $\abs{\Lambda}\log q$. The second sum can be bounded in absolute value, because by Lemma \ref{lemma:non-nullness-inequality} it holds that 
\begin{align*}
    \abs{\log \mu(\eta_{\Lambda})} \leq \abs{\Lambda}\log \frac{1}{\delta}.
\end{align*}
By taking the density limit in~\eqref{proof:attractor-property-1} we see that
\begin{align*}
    -M \leq h(\nu\lvert \mu) \leq M. 
\end{align*}
Now let $\nu^*$ be some weak limit point. If $\nu^* \notin \mathscr{G}(\gamma)$, then by the Proposition \ref{proposition:attractor-property-helper} and weak convergence there exists an open neighborhood $G$ of $\nu^*$ such that $\nu_{t_n} \in G$ for all $n \geq N(G)$ and 
\begin{align*}
    \forall \rho  \in G \ \forall 0 \leq s \leq \varepsilon: \quad h(\rho_s\lvert \mu) - h(\rho\lvert \mu) \leq -\delta s . 
\end{align*}
This implies that for all $m \in \N$
\begin{align*}
    -M \leq h(\nu^*\lvert \mu) 
    &\leq h(\nu_{t_{m+N}}\lvert \mu) - h(\nu_{t_N}\lvert \mu) + h(\nu\lvert \mu)\leq - \delta\sum_{k=0}^{m-1}\min\set{\varepsilon, t_{N+k+1}-t_{N+k}} + M. 
\end{align*}
Since $t_n$ increases to infinity, we necessarily have that the sum on the right-hand side diverges to infinity as $n$ tends to infinity. But this leads to a contradiction, since $M$ is finite. Therefore, we must have $\nu^* \in \mathscr{G}(\gamma)$. 
\end{proof}

\section{Acknowledgements}
The authors would like to thank the three anonymous referees for their insightful feedback that helped to substantially improve this manuscript. Benedikt Jahnel acknowledges the financial support of the Leibniz Association within the Leibniz Junior Research Group on \textit{Probabilistic Methods for Dynamic Communication Networks} as part of the Leibniz Competition. 

\bibliography{references}
\bibliographystyle{alpha}
\end{document}